\documentclass[fleqn,review,12pt,3p]{elsarticle}

\usepackage{amssymb}
\usepackage{graphicx} 
\PassOptionsToPackage{hyphens}{url}\usepackage{hyperref}
\usepackage{amsmath, amsfonts, amsthm}
\usepackage{booktabs, tabularx}
\usepackage{multirow}
\usepackage{tikz}
\usetikzlibrary{shapes.geometric}

\usepackage{todonotes}
\usepackage[font=footnotesize,labelfont=bf]{caption}
\usepackage[font=footnotesize,labelfont=bf]{subcaption}
\usepackage{mathtools, cuted}
\usepackage{xcolor}
\usepackage{mathtools}
\usepackage{longtable}
\allowdisplaybreaks
\usepackage{pgfplotstable}
\usepackage{pgfplots}
\usepackage{setspace}
\pgfplotsset{compat=1.18,
    axis line style={black},
    every axis label/.append style ={black},
    every tick label/.append style={black},
    title style={color=black},
    legend style={text=black},
}
\usepgfplotslibrary{fillbetween}

\definecolor{single}{RGB}{255,0,0}
\definecolor{double}{RGB}{0,0,255}

\usepackage{lineno}

\begin{document}
\begin{frontmatter}
    \title{Future-proof ship pipe routing: navigating the energy transition}  
    \author[1,2]{Berend Markhorst \corref{cor1}} \ead{berend.markhorst@cwi.nl} 
    \author[2]{Joost Berkhout} 
    \author[2]{Alessandro Zocca} 
    \author[3]{Jeroen Pruyn} 
    \author[1,2]{Rob van der Mei} 
    \cortext[cor1]{Corresponding author} 
    \affiliation[1]{Stochastics Department CWI, Science Park 123, 1098 XG, Amsterdam, The Netherlands}
    \affiliation[2]{Mathematics Department VU Amsterdam, Boelelaan 1111, 1081 HV, Amsterdam, The Netherlands}
    \affiliation[3]{Department of Maritime and Transport Technology, TU Delft, Leeghwaterstraat 2628 CD, Delft, The Netherlands}
    \date{\today}

    \begin{abstract}
    The maritime industry must prepare for the energy transition from fossil fuels to sustainable alternatives. Making ships future-proof is necessary given their long lifetime, but it is also complex because the future fuel type is uncertain. Within this uncertainty, one typically overlooks pipe routing, although it is a crucial driver for design time and costs. Therefore, we propose a mathematical approach for modeling uncertainty in pipe routing with deterministic, stochastic, and robust optimization. All three models are based on state-of-the-art integer linear optimization models for the Stochastic Steiner Forest Problem and adjusted to the maritime domain using specific constraints for pipe routing. We compare the models using both artificial and realistic instances and show that considering uncertainty using stochastic optimization and robust optimization leads to cost reductions of up to 22\% in our experiments.
\end{abstract}


\end{frontmatter}

\section{Introduction}
The maritime industry is a significant contributor to global greenhouse gas emissions as it accounts for 2-3\% of global carbon emissions \cite{international_maritime_organization_fourth_2020}. Therefore, the International Maritime Organization (IMO) and the United Nations (UN) have made regulations and guidelines for the maritime industry. To stop global warming, emissions must be reduced as soon as possible. Approximately a 40\% reduction per vessel is to be achieved by 2030 and even net zero for the fleet by 2050 \cite{international_maritime_organization_imo_2023}. 

These guidelines should stimulate the maritime industry to transition from fossil fuels to sustainable alternatives. However, ship owners optimize for the economic situation at the beginning of a ship's construction \cite{pruyn_are_2017}.  This means sustainability is often not taken into account. Yet, with a service life of 20 to 30 years at least \cite{dinu_maritime_2015}, a ship should be profitable in all (economic) conditions encountered in that period. As a result, ships are currently already sub-optimal at delivery. To deal with this issue, we must consider future alternative fuels already in the design phase to prevent this sub-optimality.

Transitioning to a new fuel type would impact the engine, the fuel storage, and the piping between the two. The first two are extensively researched already \cite{lindstad_reduction_2021, zwaginga_evaluation_2022}, and piping as well \cite{dong_ship_2022, lin_multi-objective_2023, yan_hybrid_2024}, yet it is often overlooked in the context of the energy transition. Pipe routing takes over 50\% of the total detail-design labor hours \cite{park_pipe-routing_2002}. Additionally, the labor costs yield 60\% of the total costs of a ship \cite[Section 1.2]{asmara_pipe_2013}. Hence, pipe routing greatly determines both design time and costs.

Furthermore, pipe routing constraints are heavily dependent on the corresponding fuel type \cite{lloyds_register_rules_2023}. Even though present studies fundamentally mention the same alternative fuels, the industry has not agreed on a single optimal future marine fuel \cite{prussi_potential_2021}. They state that the future mix of fuels will depend on their expected price, availability, and suitability for the specific ship.

This work focuses on the mathematical optimization of pipe routing in ship design in light of the uncertainty created by the energy transition. The goal is to find a pipe routing that minimizes the costs of installing pipes when building a ship and adjusting this routing in the future, also referred to as retrofit costs. To that end, we build upon and combine aspects of pipe routing and the mathematical optimization literature.

\paragraph{Related pipe routing literature} The existing literature on pipe routing is described in several surveys, see~\cite{xiao-long_qian_survey_2008, asmara_pipe_2013, blokland_literature_2023}. The synthesis tables from \cite[Section 5]{blokland_literature_2023} give a detailed overview of the current state-of-the-art. According to this table, the pipe routing literature mainly describes deterministic models that do not take uncertainty into account. Even more, considering the uncertainty about the energy transition can reduce future retrofit costs, which has not been done in the literature so far \cite[Section 6]{blokland_literature_2023}.

\paragraph{Related mathematical optimization literature} The goal of ship pipe routing is connecting multiple points, such as engine room(s) and tank(s), using as little material as possible while following the rules described in ~\cite{lloyds_register_rules_2023}. Mathematically, this is close to a Steiner Tree Problem (STP), which is a well-studied problem in combinatorial optimization, see~\cite{ljubic_solving_2021} for a recent overview. We represent a ship by a graph in which the vertices denote (engine) room(s) and tank(s), and the edges denote a connection between those vertices where we can install pipes. The objective is to connect a set of terminals (i.e., a subset of vertices) in a given graph using edges with minimal total costs. For example, a set of terminals may consist of an engine and multiple fuel tanks. STP is known to be NP-hard \cite[p. 208-209]{garey_computers_2009}, meaning that it is unlikely that there exists a scalable (polynomial-time) algorithm that can find the optimal solution to all STP instances.

The stochastic equivalent of the STP is called the Stochastic Steiner Tree Problem (SSTP), which considers two stages and a finite number of scenarios with the corresponding probabilities, terminal groups, and edge costs. In the first stage, it is unknown which scenario will occur in the second stage. The question is which edges to buy in the first stage and which (more expensive) edges to buy in the second stage. Approximation algorithms for the SSTP are described in~\cite{gupta_boosted_2004, hutchison_stochastic_2005, swamy_approximation_2006, charikar_stochastic_2007, gupta_lp_2007}, whereas an exact model that uses a two-stage branch-and-cut algorithm based on Benders' decomposition is discussed in~\cite{cheong_solving_2010}. For the \hyperlink{https://dimacs11.zib.de}{11th DIMACS challenge}, a genetic algorithm has been developed and is discussed in~\cite{hokama_heuristic_2014}. Additionally, a comparison between different Integer Linear Optimization (ILO) models for the SSTP is made in~\cite{zey_ilp_2016}. Finally, a two-stage branch-and-cut algorithm based on a decomposed model is discussed in~\cite{ljubic_stochastic_2017}, whereas a new decomposition model is described in~\cite{leitner_decomposition_2018}.

In this work, multiple groups of terminals (i.e., tanks and engine rooms) must be connected in ship pipe routing. Consequently, we consider a generalization of STP called the Steiner Forest Problem (SFP) \cite{kao_steiner_2008}, in which multiple terminal groups must be connected. As STP is NP-hard, also SFP is NP-hard \cite{gassner_steiner_2010}. For this problem, approximation algorithms are described in~\cite{bateni_approximation_2011, civril_approximation_2019, ghalami_approximation_2022}, whereas a greedy algorithm is discussed in~\cite{gupta_greedy_2015} and a local-search algorithm in~\cite{gros_local-search_2017}. We mainly base our work on a study \cite{schmidt_stronger_2021} that describes and compares different ILO models for the SFP. We have only found one work that studies the Stochastic Steiner Forest Problem (SSFP) with an approximation algorithm \cite{gupta_constant-factor_2009}. To the best of our knowledge, no paper in the literature describes an SSFP where the terminal groups from the first stage scenario must also be connected.

\paragraph{Our contributions} This work is the first to introduce the 2-Stage Stochastic Steiner Forest Problem (2S-SSFP), which we base on the SSTP problem formulation from~\cite{cheong_solving_2010, zey_ilp_2016, leitner_decomposition_2018}. Our work differs from these three papers as
\begin{itemize}
    \item we generalize from STP to SFP by using multiple terminal groups per scenario;
    \item we generalize the problem by requiring the terminal groups from the first stage scenario to be connected.
\end{itemize}
We study deterministic (DO), stochastic (SO), and robust optimization (RO) models adjusted with pipe routing constraints. For these models, we build upon~\cite{schmidt_stronger_2021}, as this work provides models with state-of-the-art LO relaxation bounds for SFP, leading to enhanced optimization performance in practice. To study the scalability of the models and the relative gains of considering uncertainty, we use artificial instances and a realistic instance made in collaboration with maritime experts. Although this work focuses on ship pipe routing, we would like to emphasize that the studied network design problem is also relevant in other domains such as in, for example, telecommunication \citep{cheong_solving_2010,hokama_heuristic_2014,ljubic_stochastic_2017,leitner_decomposition_2018}.

\paragraph{Outline} The rest of this paper is structured as follows. Sections~\ref{sec:deterministic_problem_formulation} and~\ref{sec:stochastic_problem_formulation} formulate the deterministic and stochastic problem, respectively, and elaborate on the mathematical models. A discussion of the experiments and the corresponding results follows in Section~\ref{sec:results}. Finally, we conclude the paper with our conclusions and directions for future research in Section~\ref{sec:conclusion}. 
\section{Deterministic problem formulation} \label{sec:deterministic_problem_formulation}
In this section, we describe the deterministic problem formulation of the SFP. This model captures the ship pipe routing problem for a single fuel type without considering the possibility of a transition to another fuel in the future, which is often the focus of current practice \cite{pruyn_are_2017, pruyn_benchmarking_2020}. We introduce this model in such a way that it can be reused in the stochastic and robust optimization models in Section~\ref{sec:stochastic_problem_formulation}. More specifically, each future scenario corresponds to an equivalent deterministic ship pipe routing problem for a particular fuel type in the future. The future scenario-dependent parameters and decision variables will use the same notation as introduced in this section but then decorated by $^{(s)}$, for example, the current admissible edges will be denoted by $E$ and in future scenario $s$ by $E^{(s)}$. 

In the following, we first give the general formulation of the SFP and explain how this general problem relates to ship pipe routing in Section~\ref{sec:descriptionSPR}. Then, we present the mathematical models of our deterministic optimization ILO model in Section~\ref{sec:deterministic_IP}.

\subsection{Problem description for ship pipe routing} \label{sec:descriptionSPR}
\begin{figure}[!h]
\centering
\begin{subfigure}{0.55\linewidth}
    \centering
    \includegraphics[width=\linewidth]{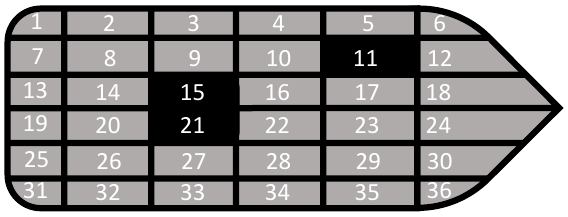}
    \caption{}
    \label{fig:1}
\end{subfigure}
\hfill
\begin{subfigure}{0.35\linewidth}
        \centering
        \resizebox{\linewidth}{!}{
    \begin{tikzpicture}[node distance={1.5cm}, mycircle/.style={draw, circle, fill=gray!20, thick, minimum size=1cm}, mydashed/.style={gray!50, dashed, line width=3}] 
        \node[mycircle] (1) {1}; 
        \node[mycircle] (2) [right of=1] {2}; 
        \node[mycircle] (3) [right of=2] {3}; 
        \node[mycircle] (4) [right of=3] {4}; 
        \node[mycircle] (5) [right of=4] {5}; 
        \node[mycircle] (6) [right of=5] {6};
        
        \node[mycircle] (7) [below of=1] {7}; 
        \node[mycircle] (8) [right of=7] {8}; 
        \node[mycircle] (9) [right of=8] {9}; 
        \node[mycircle] (10) [right of=9] {10}; 
        \node[] (11) [right of=10] {}; 
        \node[mycircle] (12) [right of=11] {12};
        
        \node[mycircle] (13) [below of=7] {13}; 
        \node[mycircle] (14) [right of=13] {14}; 
        \node[] (15) [right of=14] {}; 
        \node[mycircle] (16) [right of=15] {16}; 
        \node[mycircle] (17) [right of=16] {17}; 
        \node[mycircle] (18) [right of=17] {18};
        
        \node[mycircle] (19) [below of=13] {19}; 
        \node[mycircle] (20) [right of=19] {20}; 
        \node[] (21) [right of=20] {}; 
        \node[mycircle] (22) [right of=21] {22}; 
        \node[mycircle] (23) [right of=22] {23}; 
        \node[mycircle] (24) [right of=23] {24};
        
        \node[mycircle] (25) [below of=19] {25}; 
        \node[mycircle] (26) [right of=25] {26}; 
        \node[mycircle] (27) [right of=26] {27}; 
        \node[mycircle] (28) [right of=27] {28}; 
        \node[mycircle] (29) [right of=28] {29}; 
        \node[mycircle] (30) [right of=29] {30};
        
        \node[mycircle] (31) [below of=25] {31}; 
        \node[mycircle] (32) [right of=31] {32}; 
        \node[mycircle] (33) [right of=32] {33}; 
        \node[mycircle] (34) [right of=33] {34}; 
        \node[mycircle] (35) [right of=34] {35}; 
        \node[mycircle] (36) [right of=35] {36};
        
        \draw[mydashed] (1) -- (2);
        \draw[mydashed] (2) -- (3);
        \draw[mydashed] (3) -- (4);
        \draw[mydashed] (4) -- (5);
        \draw[mydashed] (5) -- (6);
        
        \draw[mydashed] (7) -- (8);
        \draw[mydashed] (8) -- (9);
        \draw[mydashed] (9) -- (10);
        
        \draw[mydashed] (13) -- (14);
        \draw[mydashed] (16) -- (17);
        \draw[mydashed] (17) -- (18);

        \draw[mydashed] (19) -- (20);
        \draw[mydashed] (22) -- (23);
        \draw[mydashed] (23) -- (24);
        
        \draw[mydashed] (25) -- (26);
        \draw[mydashed] (26) -- (27);
        \draw[mydashed] (27) -- (28);
        \draw[mydashed] (28) -- (29);
        \draw[mydashed] (29) -- (30);

        \draw[mydashed] (31) -- (32);
        \draw[mydashed] (32) -- (33);
        \draw[mydashed] (33) -- (34);
        \draw[mydashed] (34) -- (35);
        \draw[mydashed] (35) -- (36);

        \draw[mydashed] (1) -- (7);
        \draw[mydashed] (2) -- (8);
        \draw[mydashed] (3) -- (9);
        \draw[mydashed] (4) -- (10);
        \draw[mydashed] (6) -- (12);

        \draw[mydashed] (7) -- (13);
        \draw[mydashed] (8) -- (14);
        \draw[mydashed] (10) -- (16);
        \draw[mydashed] (12) -- (18);

        \draw[mydashed] (13) -- (19);
        \draw[mydashed] (14) -- (20);
        \draw[mydashed] (16) -- (22);
        \draw[mydashed] (17) -- (23);
        \draw[mydashed] (18) -- (24);

        \draw[mydashed] (19) -- (25);
        \draw[mydashed] (20) -- (26);
        \draw[mydashed] (22) -- (28);
        \draw[mydashed] (23) -- (29);
        \draw[mydashed] (24) -- (30);

        \draw[mydashed] (25) -- (31);
        \draw[mydashed] (26) -- (32);
        \draw[mydashed] (27) -- (33);
        \draw[mydashed] (28) -- (34);
        \draw[mydashed] (29) -- (35);
        \draw[mydashed] (30) -- (36);

    \end{tikzpicture}
}
        \caption{}
        \label{fig:2}
\end{subfigure}
\caption{Example of a ship pipe routing problem instance and the corresponding abstract representation as a graph. Figure~\ref{fig:1} shows a schematic top-down view of a ship's deck. We can install pipes in the gray rooms but not in the black rooms. Figure~\ref{fig:2} displays a graph representation of the ship's deck. Note that we omitted the three vertices whose corresponding rooms are located at the black squares of Figure~\ref{fig:1}. Dashed lines denote possible pipe connections.}
\label{fig:example_representation}
\end{figure}
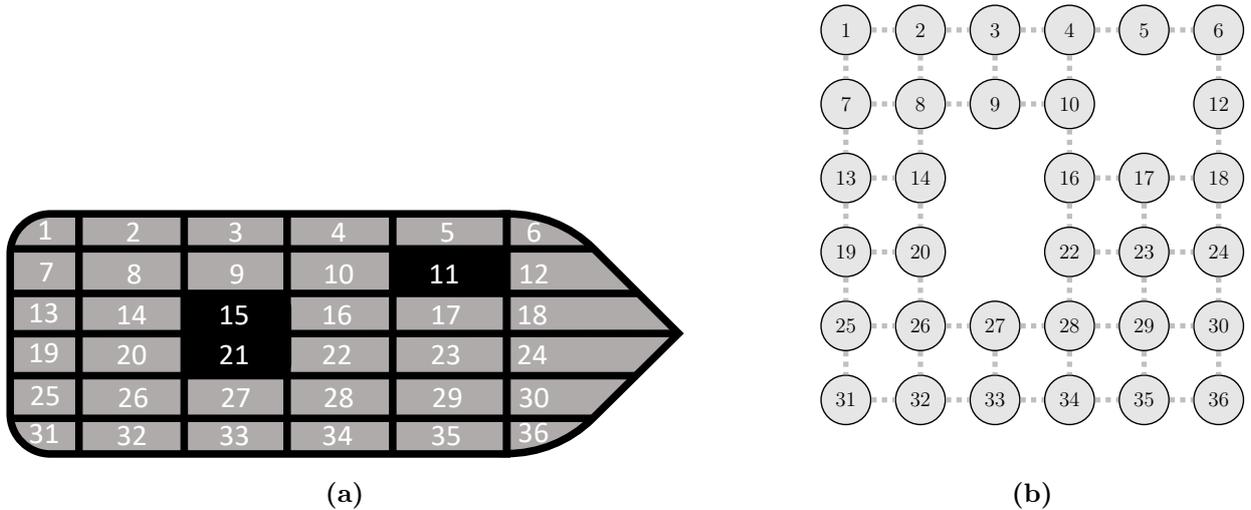

We model a ship as an undirected graph $G = (\mathcal{V}, \mathcal{E})$ where the finite vertex set $\mathcal{V}$ denotes the set of (ship) rooms and $\mathcal{E}$ represents the possible pipe connections between them, i.e., $\mathcal{E} := \left\{ (u,v): u \in \mathcal{V}, v \in \mathcal{V}, u \sim v, u < v \right\}$, where $u \sim v$ denotes adjacency between the vertices $u$ and $v$. Furthermore, $\mathcal{A}$ is the set of arcs of the bi-direction of $G$, represented by $\left\{ (u,v): u \in \mathcal{V}, v \in \mathcal{V}, u \sim v \right\}$. Each (ship) room can contain one or more engines or fuel tanks. We illustrate in Figure~\ref{fig:1} an example of a schematic top-down view of a ship's deck and in Figure~\ref{fig:2} the corresponding graph representation. Later, when we consider future scenarios, we will assume that the ship's graph representation $G$ remains the same. 

Pipes must connect each ship's engine to one or more tanks to enable fuel transportation. This information is captured by a given set of $K$ terminal groups $T = (T^k)_{k \in \mathcal{K}}$, where $T^k \subseteq {\cal V}$ denotes the $k$-th terminal group consisting of terminal vertices, in short terminals, and $\mathcal{K} = \{1, \ldots, K\}$. We must install pipes such that all vertices within a terminal group $T^k$ are connected for all $k \in \mathcal{K}$. Without loss of generality, we can take these terminal groups to be pairwise disjoint subsets of vertices, i.e., $\bigcap_{k \in \mathcal{K}} T^k = \emptyset$ \cite[Section 1]{schmidt_stronger_2021}. 

We refer to a scenario in this article as a pipe routing instance in the present or in the future (in contrast, the literature often denotes a scenario as a branch in the scenario tree, see Figure~\ref{fig:scenario_tree} for an example). We assume that every scenario (now and in the future) corresponds to one fuel type, which does not rule out the possibility that future scenarios differ from the present scenario. Different fuel types require different pipe types \cite{lloyds_register_rules_2023}. Therefore, we introduce a finite set of available pipes $\mathcal{P}$ and the subset $P \subseteq \mathcal{P}$ that describes feasible pipes for the fuel type under consideration. Introducing $\mathcal{P}$ and $P$ allows us to consider different fuel types in the future. The cost of placing a pipe $p \in \mathcal{P}$ at a particular edge $(u,v) \in \mathcal{E}$ is given $c_{puv} > 0$ and all costs are collected in $C = (c_{puv})_{p \in \mathcal{P}, (u,v) \in \mathcal{E}}$. 

Due to regulations, pipes transporting dangerous fuel types cannot cross certain rooms on a ship \cite{lloyds_register_rules_2023}. For example, diesel cannot be routed through a room adjacent to the waterside to reduce the chance of pollution in an accident. For this purpose, we introduce the set of admissible edges $E \subseteq \mathcal{E}$ for the fuel type under consideration, i.e., only edges in $E$ can be used to install pipes to connect terminal groups. Again, introducing both $E$ and $\mathcal{E}$ allows us to consider different fuel types in the future. Similarly, note that $E  \subseteq \mathcal{E}$ induces a subset of admissible arcs $A \subseteq \mathcal{A}$, i.e., only arcs in $A$ can be used to connect vertices in terminal groups. 

An instance of our problem is denoted by $I = (G, \mathcal{P}, \mathcal{T}, P, C, E)$, where $G$ (and thus $\mathcal{A}$) and $P$ are fixed, and $\mathcal{T}$, $P$, $C$, and $E$ (and thus $A$) can vary in future scenarios. We define $\mathcal{M}(I)$ as the set of all feasible solutions for instance $I$. A feasible solution $\textbf{S} \in \mathcal{M}(I)$ for a fuel type under consideration is given by a set of edge-pipe pairs that (i) connect all vertices in each terminal group from $\mathcal{T}$, (ii) uses only feasible pipes from $P$ for the terminal connections, and (iii) only uses feasible edges for the terminal connections. Note that a feasible $\textbf{S} \in \mathcal{M}(I)$ may use different feasible pipes to connect terminals and can also install infeasible pipes as long as they are not used to connect terminals (this can be efficient to anticipate for future scenarios later on). To be able to take future scenarios into account, we assume that we are given the set $\textbf{S}_0$ of edge-pipe pairs that describe which pipes already exist on which edge. Specifically, $(p, (u,v)) \in \textbf{S}_0$ means that pipe $p \in P$ is present at edge $(u,v) \in \mathcal{E}$. We let $F(I, \textbf{S}_0, \textbf{S})$ be the cost of solution $\textbf{S} \in \mathcal{M}(I)$ for instance $I$ when edge-pipe pairs $\textbf{S}_0$ are present. Then our mathematical optimization problem can be written as
\begin{equation}\label{eq:DO_problem_new}
        \textbf{\mbox{(DO)}} \qquad \min_{\textbf{S} \in \mathcal{M}(I)} F(I, \textbf{S}_0, \textbf{S}) = \min_{\textbf{S} \in \mathcal{M}(I)} \sum_{(p, (u,v)) \in \textbf{S} \setminus \textbf{S}_0}  c_{puv} .
\end{equation}
Note that our problem, with $\textbf{S}_0 = \emptyset$, reduces to a Steiner forest problem $(G, \mathcal{T})$ if it holds that $\mathcal{P} = P = \{1\}$, $\mathcal{E} = E$ (and thus $\mathcal{A} = A$) in $I = (G, \mathcal{P}, \mathcal{T}, P, C, E)$, i.e., there is only one feasible pipe type, all edges are feasible, and there are no current pipes in the graph. Consequently, (DO) is NP-hard.

\subsection{Deterministic ILO model for ship pipe routing} \label{sec:deterministic_IP}
We can reformulate (DO) to an integer linear optimization (ILO) model so that (commercial) ILO solvers can solve our problem. This section introduces the deterministic ILO model and the corresponding decision variables and constraints. For the deterministic ILO model, we build on the work of~\cite{schmidt_stronger_2021}, as this work describes the cut-based and flow-based ILO models for the SFP, which are equivalent. In this work, we use both the undirected and directed flow-based ILO models from~\cite{schmidt_stronger_2021} tailored to the ship pipe routing problem. An overview of the ILO models can be found in Table~\ref{tab:overview_labels_of_our_models}, which shows that we consider three optimization types (DO, RO, and SO) and two model types regarding the flow (undirected and directed). For example, (DO-U) and (DO-D) are ILO models for DO using undirected and directed formulations, respectively.  We explain later the difference between the two model types.

\begin{table}[ht]
    \centering
    \caption{Abbreviations of our different ILO models. Model types specify whether directed flows or undirected flows are considered.}
    \label{tab:overview_labels_of_our_models}
    \begin{tabular}{@{}lll@{}}
        \toprule
        Optimization type \textbackslash \; model type & Undirected (U) & Directed (D) \\\midrule
        Deterministic optimization (DO) & (DO-U) & (DO-D) \\
        Robust optimization (RO) & (RO-U) & (RO-D) \\
        Stochastic optimization (SO) & (SO-U) & (SO-D) \\
        \bottomrule
    \end{tabular}
\end{table}

\paragraph{Deterministic optimization with an undirected formulation}
For (DO-U) in~\eqref{eq:basic_deterministic}, we introduce a binary decision variable $x_{puv}$ that equals $1$ if pipe $p \in \mathcal{P}$ is installed at edge $(u,v) \in \mathcal{E}$, and $0$ else. To ensure that pipes connect all terminals in each terminal group, we will let an artificial flow go through the pipes. In particular, in each terminal group $k \in \mathcal{K}$, we randomly (without losing on generality) designate one terminal as the root vertex $r_k \in T^k$ of that particular terminal group for orientation purposes. The set of all root vertices is defined as $\mathcal{R} = \{r_1, \ldots r_K\}$. From the root vertex, we let an (artificial) flow go to all remaining terminals in the corresponding terminal group. We define $f_{tpuv}$ as the flow amount to terminal $t \in \mathcal{T} \setminus \mathcal{R}$ from the corresponding root over pipe $p \in P$ at arc $(u,v) \in A$. The formulation of (DO-U) is called ``undirected'' as it does not force a flow direction upon the installed pipes for different terminals. Table~\ref{tab:notation_overview1a} gives the comprehensive list of the sets, parameters, and decision variables used to formulate (DO-U) in~\eqref{eq:basic_deterministic}.

\begin{longtable}{@{}p{0.2\textwidth}p{0.8\textwidth}@{}} 
\caption{Notation overview for our ship pipe routing problem (DO) and its ILO model (DO-U).}\label{tab:notation_overview1a}\\
\toprule
Sets &  \\ 
\midrule
$\mathcal{V}$ & Set of vertices.\\
$\mathcal{E}$ & Set of edges, represented by $\{(u,v): u \in \mathcal{V}, v \in \mathcal{V}, u \sim v, u < v\}$.\\
$E$ & Set of admissible edges; $E \subseteq \mathcal{E}$.\\
$\mathcal{A}$ & Set of arcs, represented by $\left\{ (u,v): u \in \mathcal{V}, v \in \mathcal{V}, u \sim v \right\}$.\\
$A$ & Set of admissible arcs; $A \subseteq \mathcal{A}$.\\
$\mathcal{K}$ & Set of indices for the terminal groups, indexed from $1$ to $K$.\\
$\mathcal{T}$ & Set of all terminal groups $\mathcal{T} = (T^k)_{k \in \mathcal{K}}$, where $T^k \subseteq \mathcal{V}$.\\
$\mathcal{R}$ & Set of root vertices; $\mathcal{R} = \{r^1, \ldots r^k\}$, where $r^k \in T^k$ for terminal group $k \in \mathcal{K}$. Note that the root vertex is chosen arbitrarily for each terminal group.\\
$\mathcal{P}$ & Set of available pipe types.\\
$P$ & Set of feasible pipe types; $P \subseteq \mathcal{P}$.\\
\midrule
Parameters & \\
\midrule
$c_{puv}$ & Cost parameter for installing pipe type $p \in \mathcal{P}$ at edge $(u,v) \in \mathcal{E}$.\\
\midrule
Other notation & \\
\midrule
$\tau(t)$ & The index of the unique terminal group that contains the terminal $t \in \mathcal{T}$.\\
\midrule
\multicolumn{2}{@{}l@{}}{Decision Variables} \\
\midrule
$x_{puv}$ & Binary variable equal to $1$ if we install a pipe of type $p \in \mathcal{P}$ on edge $(u,v) \in \mathcal{E}$, and $0$ otherwise.\\
$f_{tpuv}$ & Binary variable equal to $1$ if there is a flow over arc $(u,v) \in A$ with pipe type $p \in P$ for a route between vertex $t \in \mathcal{T} \setminus R$ and $r^{\tau(t)}$, and $0$ otherwise.\\
\bottomrule
\end{longtable}

\begin{subequations}\label{eq:basic_deterministic}
    \begin{alignat}{3}
        \mbox{\textbf{(DO-U)}} && \notag\\
        \min \quad & \sum_{(p,(u,v)) \in (\mathcal{P} \times \mathcal{E}) \setminus \textbf{S}_0} x_{puv} \cdot c_{puv} & \label{eq:basic_do1} \\
        \begin{split}
            \mbox{s.t.} \quad & \sum_{p \in P} \left (\sum_{u: (v, u) \in A} f_{tpvu} - \right.\\
            & \left. \sum_{u: (u, v) \in A} f_{tpuv} \right ) = \left\{\begin{matrix*}[l]
                1 & \text{if } v = r^{\tau(t)} \\
                -1 & \text{if } v = t \\
                0 &  \text{otherwise}\\
                \end{matrix*}
            \right. 
        \end{split}
        & \qquad & \forall v \in \mathcal{V}, \forall t \in \mathcal{T} \setminus \mathcal{R} \label{eq:basic_do2} \\
        & f_{tpuv} + f_{tpvu} \leq x_{puv} & \qquad & \forall t \in \mathcal{T} \setminus \mathcal{R}, \forall p \in P, \forall (u,v) \in E \label{eq:basic_do3}\\
        & f_{tpuv} \in \mathbb{B} & \qquad & \forall t \in \mathcal{T} \setminus \mathcal{R}, \forall p \in P, \forall (u,v) \in A \label{eq:basic_do4}\\
        & x_{puv} \in \mathbb{B} & \qquad & \forall p \in \mathcal{P}, \forall (u,v) \in \mathcal{E} \label{eq:basic_do5}
    \end{alignat}
\end{subequations}

The goal of~\eqref{eq:basic_deterministic} is to minimize the costs of connecting the terminals within all terminal groups~\eqref{eq:basic_do1} while ensuring the conservation of flows as formulated in~\eqref{eq:basic_do2} using only admissible edges and feasible pipe types. Note that we added the summation over the set of pipes to allow for different (feasible) pipe types within a connection between a root and terminal. We ensure that an edge is used only in one direction for each connection between $r^{\tau(t)}$ and terminal $t$ and connect the decision variables $f_{tpuv}$ and $x_{puv}$ in~\eqref{eq:basic_do3}. Note that we restrict to feasible pipe types and admissible edges in our flow decision variable $f_{tpuv}$ and that~\eqref{eq:basic_do5} can be relaxed because~\eqref{eq:basic_do4} enforces the integrality of $x_{puv}$. For simplicity, and with a slight abuse of notation, we may capture all feasible $\textbf{x}$ and $\textbf{f}$ for (DO-U) in $\mathcal{M}(I)^{\mbox{(DO-U)}}$, so that (DO-U) can be compactly written as
\begin{subequations} \label{eq:basic_deterministic_short}
    \begin{alignat}{3}
        \mbox{\textbf{(DO-U)}} \qquad & \min \quad & \sum_{(p,(u,v)) \in (\mathcal{P} \times \mathcal{E}) \setminus \textbf{S}_0} x_{puv} \cdot c_{puv} & \label{eq:basic_do1b} \\
        & \mbox{s.t.} \quad & \textbf{x}, \textbf{f} \in \mathcal{M}(I)^{\mbox{(DO-U)}} & \label{eq:basic_do2b} 
    \end{alignat}
\end{subequations}

\begin{figure}[ht]
    \centering
    \begin{minipage}[ht]{0.30\linewidth}
    \centering
    \begin{tikzpicture}[node distance={15mm}, thick, mycircle/.style={draw, circle}] 
    \node[mycircle, color = red, fill=red] (1) {}; 
    \node[mycircle, color=blue, fill = blue] (2) [right of=1] {};
    \node[mycircle, color=red] (3) [below of=2] {};
    \node[mycircle, draw=blue] (4) [left of=3] {};
    \draw[color=black, line width=1, dashed] (1) -- (2);
    \draw[color=black, line width=1, dashed] (2) -- (3);
    \draw[color=black, line width=1, dashed] (3) -- (4);
    \draw[color=black, line width=1, dashed] (1) -- (4);
    \end{tikzpicture}
    \subcaption{}
    \label{fig:demo_instance}
    \end{minipage}
    \hfill
    \begin{minipage}[ht]{0.30\linewidth}
    \centering
    \begin{tikzpicture}[node distance={15mm}, thick, mycircle/.style={draw, circle}] 
    \node[mycircle, color = red, fill=red] (1) {}; 
    \node[mycircle, color=blue, fill = blue] (2) [right of=1] {};
    \node[mycircle, color=red] (3) [below of=2] {};
    \node[mycircle, draw=blue] (4) [left of=3] {};

    \draw[color=black, line width=1] (1) -- (2);
    \draw[color=black, line width=1] (2) -- (3);
    \draw[color=black, line width=1] (3) -- (4);
    \draw[color=black, line width=1] (1) -- (4);
    
    \draw[color=red, line width=2, style=->] (1) to[out=90,in=90] (2);
    \draw[color=red, line width=2, style=->] (2) to[out=0,in=0] (3);
    \draw[color=red, line width=2, style=->] (4) to[out=-90,in=-90] (3);
    \draw[color=red, line width=2, style=->] (1) to[out=180,in=180] (4);

    \draw[color=blue, line width=2, style=->] (2) to[out=-90,in=-90] (1);
    \draw[color=blue, line width=2, style=->] (1) to[out=0,in=0] (4);
    \draw[color=blue, line width=2, style=->] (3) to[out=90,in=90] (4);
    \draw[color=blue, line width=2, style=->] (2) to[out=180,in=180] (3);

    \end{tikzpicture}
    \subcaption{}
    \label{fig:basic_demo}
    \end{minipage}
    \hfill
    \begin{minipage}[ht]{0.30\linewidth}
    \centering
    \begin{tikzpicture}[node distance={15mm}, thick, mycircle/.style={draw, circle}] 
    \node[mycircle, color = red, fill=red] (1) {}; 
    \node[mycircle, color=blue, fill = blue] (2) [right of=1] {};
    \node[mycircle, color=red] (3) [below of=2] {};
    \node[mycircle, draw=blue] (4) [left of=3] {};
    \draw[color=red, line width=2, style=->] (1) to[out=90,in=90] (2);
    \draw[color=red, line width=2, style=->] (2) to[out=0,in=0] (3);
    \draw[color=red, line width=2, style=->] (1) to[out=180,in=180] (4);

    \draw[color=black, line width=1] (1) -- (2);
    \draw[color=black, line width=1] (2) -- (3);
    \draw[color=black, line width=1] (1) -- (4);
    
    \end{tikzpicture}
    \subcaption{}
    \label{fig:advanced_demo}
    \end{minipage}
    \caption{Example that shows the difference between the (DO-U) and (DO-D). Figure~\ref{fig:demo_instance} shows a small example from Figure 1 and 2 of~\cite{schmidt_stronger_2021} with two terminal groups, red and blue, and possible connections between vertices denoted by dashed lines. The roots are denoted by circles that are filled with colors. The remaining vertices are the terminals. The costs of using an edge equal $1$ for all edges in the graph. The optimal integer solution yields an objective of $3$ because it requires three edges to connect all terminals. Figure~\ref{fig:basic_demo} shows a feasible fractional solution for (DO-U) where the colored arcs denote flows of $0.5$ for the corresponding root-terminal pairs, and the corresponding $\textbf{x}$ decision variables equal $0.5$ and are represented by the solid black lines. This solution is infeasible for (DO-D) because there are directed cycles between the two roots and the two terminals. Figure~\ref{fig:advanced_demo} shows a feasible solution for (DO-D) where the $\textbf{x}$ and $\textbf{f}$ decision variables both equal $1$. The red arcs represent the arborescence that connects both terminal groups, which yields a solution with three edges.}
    \label{fig:demo}
\end{figure}
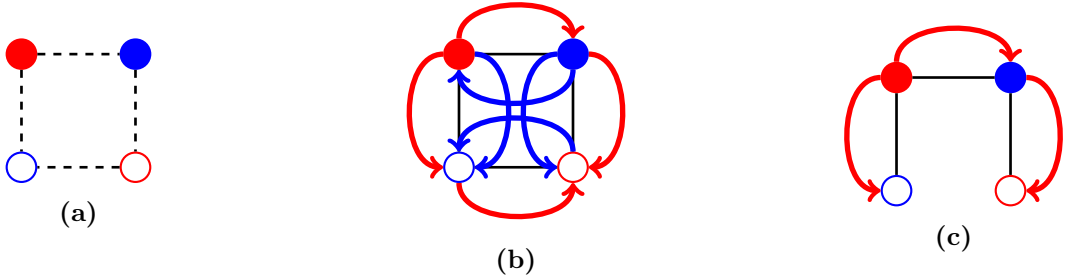

\paragraph{Deterministic optimization with a directed formulation} 
To solve an ILO model efficiently, it is important to have an ILO model that gives sharp LO-relaxation bounds \cite{schmidt_stronger_2021}. The formulation of (DO-U) can be sharpened to give better LO-relaxation bounds. The LO-relaxation of (DO-U) now allows for directed cycles of flows of different terminals, which is not tight. These directed cycles can appear for different terminals from one terminal group or when two (or more) terminal groups are connected/overlapping in the solution. By ruling out these directed cycles of flow in the LO-relaxation, a stronger model is obtained. See Figure~\ref{fig:basic_demo} for an illustrative example. 

To rule out the directed cycles of flow, we introduce an ILO model (DO-D) that ``directs'' the flow into consistent orientations. More specifically, (DO-D) improves over (DO-U) by 1) dynamically combining overlapping terminal groups effectively into one large terminal group and 2) finding an arborescence (a directed tree) for this combined terminal group in which one root is designated to send (artificial) flow to all other terminals in the combined terminal group. In (DO-D), it is ensured that all flow in an arborescence has the same orientation in the LO-relaxation. As a result, directed cycles are eliminated, as shown in Figure~\ref{fig:advanced_demo}. Hence, the LO-relaxation from (DO-D) will be tighter than the LO-relaxation from (DO-U). More details can be found in~\cite{schmidt_stronger_2021}.

For (DO-D), the decision variable $x_{puv}$ from (DO-U) remains the same. Compared to (DO-U), the flow decision variable is extended by a terminal group index $k \in \mathcal{K}$ and the resulting $f_{ktpuv}$ is the flow sent from the root of terminal group $k$ to terminal $t$ via pipe $p$ at arc $(u,v)$. Binary decision variable $z_{kl}$ is 1 when the root of the terminal group $k$ sends flow to all terminals of the terminal group $l$, and 0 else. When $z_{kl} = 1$, binary decision variable $y_{kpuv}$ equals 1 when flow from the root of terminal group $k$ is sent over arc $(u,v)$ through pipe $p$, and 0 else. Lastly, binary decision variable $y_{puv}$ equals 1 when pipe $p$ at arc $(u,v)$ is used to send flow over by any of the created arborescences, and 0 else. An overview of all decision variables is given in Table~\ref{tab:notation_overview1b}, which builds upon the notation from Table~\ref{tab:notation_overview1a}.

\begin{longtable}{@{}p{0.2\textwidth}p{0.8\textwidth}@{}} 
\caption{Notation overview for the (DO-D) problem formulation and ILO model, which builds on the notation from Table~\ref{tab:notation_overview1a}.}\label{tab:notation_overview1b}\\
\toprule
Sets &  \\ 
\midrule
$\mathcal{Q}$ & Set of non-terminal vertices, also called Steiner vertices: $\mathcal{Q} = \mathcal{V} \setminus \mathcal{T}$.\\
$\mathcal{T}^{k\ldots K}$ & Set of some terminal groups: $\mathcal{T}^{k\ldots K} = (T^k)_{k \in \{k, \ldots, \mathcal{K}\}}$.\\
$\mathcal{T}^{k\ldots K}_r$ & Set of some terminal groups without the corresponding root vertex: $\mathcal{T}^{k\ldots K}_r = \mathcal{T}^{k\ldots K} \setminus \{r^k\}$.\\
\midrule
\multicolumn{2}{@{}l@{}}{Decision Variables} \\
\midrule
$f_{ktpuv}$ & Binary variable equal to $1$ if there is a flow over arc $(u,v) \in A$ with pipe type $p \in P$ for a route between root $r^k$ of terminal group $k \in \mathcal{K}$ and terminal $t \in \mathcal{T}^{k\ldots K}_r$, and $0$ otherwise.\\
$y_{kpuv}$ & Binary variable equal to $1$ if there is a flow over arc $(u,v) \in A$ with pipe of type $p \in P$ from root $r^k$ of terminal group $k \in \mathcal{K}$, and $0$ otherwise.\\
$y_{puv}$ & Binary variable equal to $1$ if there is a flow over arc $(u,v) \in A$ with pipe of type $p \in P$ , and $0$ otherwise.\\
$z_{kl}$ & Binary variable equal to $1$ if $T^k$ and $T^l$ are in the same arborescence and $r^k$ is the root for both terminal groups, and $0$ otherwise.\\
\bottomrule
\end{longtable}

Using the extra notation, we get the following ILO model (DO-D):

\begin{subequations} \label{eq:advanced_deterministic}
    \begin{alignat}{3}
        \mbox{\textbf{(DO-D)}}&& \notag\\
        \min \quad & \sum_{(p,(u,v)) \in (\mathcal{P} \times \mathcal{E}) \setminus \textbf{S}_0} x_{puv} \cdot c_{puv} \label{eq:advanced_do1}\\
        \begin{split}
        \mbox{s.t.} \quad & \sum_{p \in P} \Big ( \sum_{u: (v, u) \in A} f_{ktpvu} \\
        & - \sum_{u: (u, v) \in A} f_{ktpuv} \Big ) =
        \left\{\begin{matrix*}[l]
            z_{kl} & \text{if } v = r^{k} \\
            -z_{kl} & \text{if } v = t \\
            0 &  \text{otherwise}
            \end{matrix*}
        \right.
        \end{split} & \qquad & 
        \left \{ \begin{matrix*}[l]
            \forall k \in \mathcal{K}, \forall t \in \mathcal{T}^{k\ldots K}_r\\
            \forall v \in \mathcal{V} \text{ with } \tau(t) = l
        \end{matrix*}
        \right. \label{eq:advanced_do2}\\
        & f_{ktpuv} \leq y_{kpuv} & \qquad & 
        \left \{ 
            \begin{matrix*}[l]
                    \forall k \in \mathcal{K}, \forall t \in \mathcal{T}^{k\ldots K}_r\\
                    \forall p \in P, \forall (u,v) \in A 
            \end{matrix*}
        \right. \label{eq:advanced_do3}\\
        & \sum_{k \in \mathcal{K}} y_{kpuv} \leq y_{puv} & \qquad & \forall p \in P, \forall (u,v) \in A \label{eq:advanced_do4}\\
        & y_{puv} + y_{pvu} \leq x_{puv} & \qquad & \forall p \in P, \forall (u,v) \in E \label{eq:advanced_do7}\\
        & \sum_{l=1}^k z_{lk} = 1 & \qquad & 
            \forall k \in \mathcal{K} \label{eq:advanced_do10}\\
        & z_{kk} \geq z_{kl} & \qquad & 
        \left \{
            \begin{matrix*}[l]
                \forall k \in \mathcal{K} \setminus \{1,K\}\\
                \forall l \in \mathcal{K} \text{ if } l \geq k + 1 
            \end{matrix*}
        \right. \label{eq:advanced_do11}\\
        & \sum_{p \in P} \sum_{u: (u, v) \in A} y_{puv} \leq 1& \qquad & \forall v \in \mathcal{V} \label{eq:advanced_do5}\\
        & \sum_{p \in P} \sum_{u: (u, t) \in A} y_{kpuv} = 0 & \qquad & \forall k \in \mathcal{K} \setminus \{1\}, \forall t \in \mathcal{T}^{1 \ldots k-1} \label{eq:advanced_do6}\\
        & \sum_{p \in P} \sum_{u: (t, u) \in A} f_{ktpuv} = 0 & \qquad & \forall k \in \mathcal{K}, \forall t \in \mathcal{T}^{k\ldots K}_r \label{eq:advanced_do16}\\
        & \sum_{p \in P} \sum_{u: (u, v) \in A} y_{puv} \leq \sum_{p \in P} \sum_{u: (v, u) \in \mathcal{A}} y_{puv} & \qquad & \forall v \in \mathcal{Q} \label{eq:advanced_do19}\\
        & \sum_{p \in P} \sum_{u: (u, v) \in A} y_{kpuv} \leq \sum_{p \in P} \sum_{u: (v, u) \in \mathcal{A}} y_{kpuv} & \qquad & \forall k \in \mathcal{K}, \forall v \in \mathcal{V} \setminus \mathcal{T}^{k\ldots K}_r \label{eq:advanced_do17}\\
        & \sum_{u: (u, r^l) \in A} y_{kpu r^l} \leq z_{kl} & \qquad &
        \left \{ \begin{matrix*}[l]
            \forall k \in \mathcal{K} \setminus {K}\\
            \forall l \in \mathcal{K} \text{ if } l \geq k + 1\\
            \forall p \in P.
        \end{matrix*}
        \right. \label{eq:advanced_do18}\\
        & f_{ktpuv} \in \mathbb{B} & \qquad & 
        \left \{
            \begin{matrix*}
                \forall k \in \mathcal{K}, \forall t \in \mathcal{T}^{k\ldots K}_r\\
                \forall p \in P, \forall (u,v) \in A                
            \end{matrix*}
        \right. \label{eq:advanced_do12}\\
        & x_{puv} \in \mathbb{B} & \qquad & \forall p \in \mathcal{P}, \forall (u,v) \in \mathcal{E} \label{eq:advanced_do20}\\
        & y_{puv} \in \mathbb{B} & \qquad & \forall p \in P, \forall (u,v) \in A \label{eq:advanced_do13}\\
        & y_{kpuv} \in \mathbb{B} & \qquad & 
        \left \{
            \begin{matrix*}
                \forall k \in \mathcal{K}, \forall p \in P\\
                \forall (u,v) \in A
            \end{matrix*}
        \right. \label{eq:advanced_do14}\\
        & z_{kl} \in \mathbb{B} & \qquad & \forall k \in \mathcal{K}, \forall j \in \{k \ldots K\}\label{eq:advanced_do15}
    \end{alignat}
\end{subequations}

The objective~\eqref{eq:advanced_do1} is similar to~\eqref{eq:basic_do1}, whereas the conservation of flows in~\eqref{eq:advanced_do2} differs from~\eqref{eq:basic_do2}. This constraint ensures that each terminal is contained in an arborescence rooted at $r^k$ for some $k \in \mathcal{K}$. From each arborescence root $r^k$, an artificial flow is sent to all remaining terminals in the arborescence. Decision variables $f_{ktpuv}$ activate $y_{kpuv}$ in~\eqref{eq:advanced_do3} whenever flow is sent from root $r^k$ to terminal $t$ via pipe $p$ at arc $(u,v)$. We ensure that every arc is part of at most one arborescence in~\eqref{eq:advanced_do4}. In case of an overlap, the corresponding arborescences are forced to be merged into one arborescence. Similarly to~\eqref{eq:basic_do3}, \eqref{eq:advanced_do7} allows for only one direction on an edge. We enforce that every terminal group is rooted at exactly one root in~\eqref{eq:advanced_do10}, whereas~\eqref{eq:advanced_do11} enforces exactly one root in each arborescence. 

Constraints~\eqref{eq:advanced_do5}-\eqref{eq:advanced_do18} are not necessary for (DO-D) to produce feasible solutions but are introduced to strengthen the model's LO-relaxation \cite{schmidt_stronger_2021}.
In~\eqref{eq:advanced_do5}, we ensure that every vertex receives flow over at most one pipe. According to the definition of $z_{kl}$, the root $r^k$ can be only responsible for terminal groups $l \geq k$. Consequently, \eqref{eq:advanced_do6} prevents a connection between root $r^k$ and a terminal from $\mathcal{T}^{1 \ldots k-1}$. Constraints~\eqref{eq:advanced_do16} prevent a flow from leaving a terminal. We denote flow-balance constraints in~\eqref{eq:advanced_do19} and~\eqref{eq:advanced_do17}, also mentioned in~\cite[Section 2.2]{leitner_decomposition_2018} for the SSTP, which state that the in-degree of a Steiner vertex cannot be larger than its out-degree: \eqref{eq:advanced_do19} enforces this in the overall solution whereas~\eqref{eq:advanced_do17} focuses on each terminal group. We enforce that the arborescence rooted at $r^k$ can use root $r^l$ if and only if $z_{kl}=1$ in~\eqref{eq:advanced_do18}. 

Finally, we include integrality constraints in~\eqref{eq:advanced_do12}-\eqref{eq:advanced_do15}. As~\eqref{eq:advanced_do12} already enforces integrality on $y_{kpuv}$, constraints \eqref{eq:advanced_do13} and~\eqref{eq:advanced_do14} can be relaxed. 

For simplicity, and with a slight abuse of notation, we may capture all feasible $\textbf{x}$, $\textbf{f}$, $\textbf{y}$, and $\textbf{z}$ for (DO-D) in $\mathcal{M}(I)^{\mbox{(DO-D)}}$, so that (DO-D) can be written as
\begin{subequations} \label{eq:advanced_deterministic_short}
    \begin{alignat}{3}
        \mbox{\textbf{(DO-D)}} \qquad & \min \quad & \sum_{(p,(u,v)) \in (\mathcal{P} \times \mathcal{E}) \setminus \textbf{S}_0} x_{puv} \cdot c_{puv} & \label{eq:advanced_do1b} \\
        & \mbox{s.t.} \quad & \textbf{x}, \textbf{f}, \textbf{y}, \textbf{z} \in \mathcal{M}(I)^{\mbox{(DO-D)}} & \label{eq:advanced_do2b} 
    \end{alignat}
\end{subequations}

\section{Accounting for uncertainty: two new problem formulations} \label{sec:stochastic_problem_formulation}
In this section, we explain how uncertainty affects the ship pipe routing problem and introduce two new optimization models to deal with it. Stochastic optimization (SO) and robust optimization (RO) are techniques that address uncertainties and variability in real-world optimization problems. They are different approaches for dealing with uncertainty; SO requires more detailed distributional information and focuses on the average case, whereas RO requires information on the support of the uncertain parameters and focuses on the worst case. Although both models are well known in the mathematical literature, applying these models in our practical context is novel \cite[Section 6]{blokland_literature_2023}. For more details about SO, also known as Stochastic Programming in the literature, we refer to~\cite{birge_introduction_2011, klein_haneveld_stochastic_2020}, whereas \cite{ben-tal_robust_2009, gorissen_practical_2015} provide more information on RO.

\subsection{Benefits of using SO and RO for pipe routing} 
\begin{figure}[ht]
    \centering
    \resizebox{\linewidth}{!}{
        \begin{tikzpicture}[node distance={35mm}, thick] 
        \node[color=black] (1) {Diesel};
        \node[color=black] (2) [below left of=1] {Methanol: $s=2$};
        \node[color=black] (3) [left of=2] {Diesel: $s=1$};
        \node[color=black] (4) [below right of=1] {$\hdots$};
        \node[color=black] (5) [right of=4] {LNG: $s=S$};
        \node[color=black] (6) [left of=3] {\textbf{Future: Second stage}};
        \node[color=black] (7) [above left of=2] {\textbf{Present: First stage}};
        \draw[color=black] (1) -> (2);
        \draw[color=black] (1) -> (3);
        \draw[color=black] (1) -> (4);
        \draw[color=black] (1) -> (5);
        \end{tikzpicture}    
    }
    \caption{Two-stage scenario tree that schematically represents the problem setting we are studying.}
    \label{fig:scenario_tree}
\end{figure}
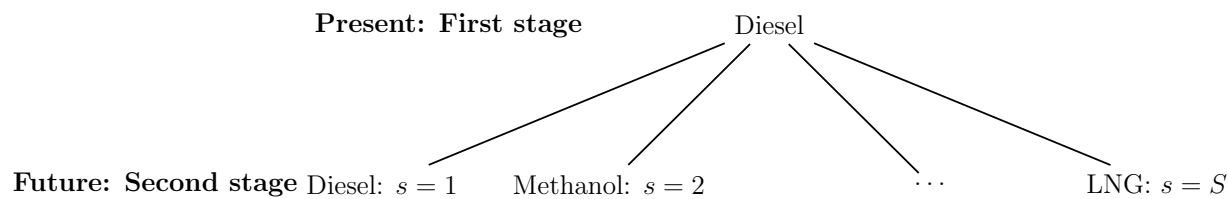

Currently, diesel is the most used ship fuel \cite{prussi_potential_2021}. As part of the energy transition, guidelines from the IMO \cite{international_maritime_organization_imo_2023} stimulate the maritime industry to transition to alternative, less polluting fuels. However, no single alternative fuel is recommended to be used in the future. The future mix of fuels used depends on many external factors, such as technology improvements, availability, and future costs \cite{prussi_potential_2021}. In our framework, we consider two periods in time, the present and the future, and refer to them as the first stage and second stage, respectively. Figure~\ref{fig:scenario_tree} shows the scenario tree schematically representing our problem setting. At each stage, pipe routes can be changed by installing extra pipes. We assume that it is not necessary to remove unused pipes. In the first stage, we consider one scenario, typically diesel. However, in the second stage, we consider multiple scenarios, each corresponding to a different (future) fuel type and possibly different pipe routes. For example, a future scenario could be methanol, which, unlike diesel, can be routed through rooms next to the waterside and requires double-walled pipes \cite{lloyds_register_rules_2023}. Because of these different characteristics, methanol might need different pipes and routes than diesel.

Ship designers optimize for the economic situation of a ship at the start of construction \cite{pruyn_are_2017, pruyn_benchmarking_2020}. As a result, ships are already sub-optimal at delivery due to ongoing technological developments. Changing the ship's pipe network after construction is not preferable. For example, retrofitting a ship from diesel to methanol is expensive because of the ship's downtime and the required complex maintenance. To make ships future-proof for both economic and sustainability reasons, alternative fuels must be taken into account in the design phase. 

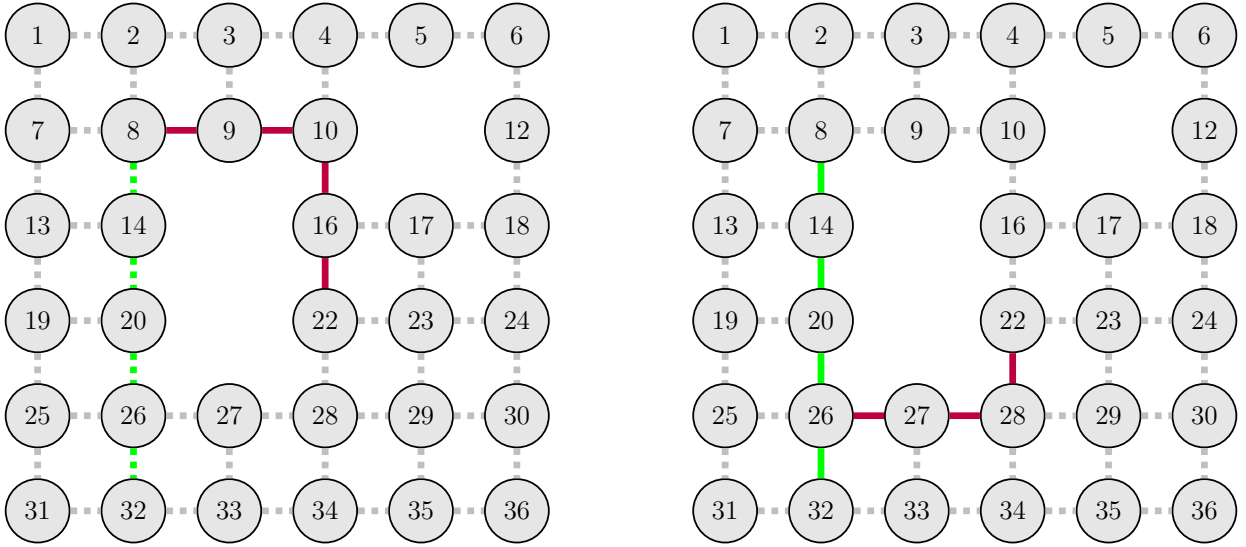
\begin{figure}[ht]
\centering
\begin{subfigure}{0.45\linewidth}
    \centering
    \resizebox{\linewidth}{!}{
    \begin{tikzpicture}[node distance={1.5cm}, mycircle/.style={draw, circle, fill=gray!20, thick, minimum size=1cm}, mydashed/.style={gray!50, dashed, line width=3}, mydashed2/.style={gray, dashed, line width=3, black}, mysolid/.style={black, line width=3}, singlewalled1/.style={purple, line width=3}, doublewalled1/.style={green, line width=3}, doublewalled2/.style={green, dashed, line width=3}] 
        \node[mycircle] (1) {1}; 
        \node[mycircle] (2) [right of=1] {2}; 
        \node[mycircle] (3) [right of=2] {3}; 
        \node[mycircle] (4) [right of=3] {4}; 
        \node[mycircle] (5) [right of=4] {5}; 
        \node[mycircle] (6) [right of=5] {6};
        
        \node[mycircle] (7) [below of=1] {7}; 
        \node[mycircle] (8) [right of=7] {8}; 
        \node[mycircle] (9) [right of=8] {9}; 
        \node[mycircle] (10) [right of=9] {10}; 
        \node[] (11) [right of=10] {}; 
        \node[mycircle] (12) [right of=11] {12};
        
        \node[mycircle] (13) [below of=7] {13}; 
        \node[mycircle] (14) [right of=13] {14}; 
        \node[] (15) [right of=14] {}; 
        \node[mycircle] (16) [right of=15] {16}; 
        \node[mycircle] (17) [right of=16] {17}; 
        \node[mycircle] (18) [right of=17] {18};
        
        \node[mycircle] (19) [below of=13] {19}; 
        \node[mycircle] (20) [right of=19] {20}; 
        \node[] (21) [right of=20] {}; 
        \node[mycircle] (22) [right of=21] {22}; 
        \node[mycircle] (23) [right of=22] {23}; 
        \node[mycircle] (24) [right of=23] {24};
        
        \node[mycircle] (25) [below of=19] {25}; 
        \node[mycircle] (26) [right of=25] {26}; 
        \node[mycircle] (27) [right of=26] {27}; 
        \node[mycircle] (28) [right of=27] {28}; 
        \node[mycircle] (29) [right of=28] {29}; 
        \node[mycircle] (30) [right of=29] {30};
        
        \node[mycircle] (31) [below of=25] {31}; 
        \node[mycircle] (32) [right of=31] {32}; 
        \node[mycircle] (33) [right of=32] {33}; 
        \node[mycircle] (34) [right of=33] {34}; 
        \node[mycircle] (35) [right of=34] {35}; 
        \node[mycircle] (36) [right of=35] {36};
        
        \draw[mydashed] (1) -- (2);
        \draw[mydashed] (2) -- (3);
        \draw[mydashed] (3) -- (4);
        \draw[mydashed] (4) -- (5);
        \draw[mydashed] (5) -- (6);
        
        \draw[mydashed] (7) -- (8);
        \draw[singlewalled1] (8) -- (9);
        \draw[singlewalled1] (9) -- (10);
        
        \draw[mydashed] (13) -- (14);
        \draw[mydashed] (16) -- (17);
        \draw[mydashed] (17) -- (18);

        \draw[mydashed] (19) -- (20);
        \draw[mydashed] (22) -- (23);
        \draw[mydashed] (23) -- (24);
        
        \draw[mydashed] (25) -- (26);
        \draw[mydashed] (26) -- (27);
        \draw[mydashed] (27) -- (28);
        \draw[mydashed] (28) -- (29);
        \draw[mydashed] (29) -- (30);

        \draw[mydashed] (31) -- (32);
        \draw[mydashed] (32) -- (33);
        \draw[mydashed] (33) -- (34);
        \draw[mydashed] (34) -- (35);
        \draw[mydashed] (35) -- (36);

        \draw[mydashed] (1) -- (7);
        \draw[mydashed] (2) -- (8);
        \draw[mydashed] (3) -- (9);
        \draw[mydashed] (4) -- (10);
        \draw[mydashed] (6) -- (12);

        \draw[mydashed] (7) -- (13);
        \draw[doublewalled2] (8) -- (14);
        \draw[singlewalled1] (10) -- (16);
        \draw[mydashed] (12) -- (18);

        \draw[mydashed] (13) -- (19);
        \draw[doublewalled2] (14) -- (20);
        \draw[singlewalled1] (16) -- (22);
        \draw[mydashed] (17) -- (23);
        \draw[mydashed] (18) -- (24);

        \draw[mydashed] (19) -- (25);
        \draw[doublewalled2] (20) -- (26);
        \draw[mydashed] (22) -- (28);
        \draw[mydashed] (23) -- (29);
        \draw[mydashed] (24) -- (30);

        \draw[mydashed] (25) -- (31);
        \draw[doublewalled2] (26) -- (32);
        \draw[mydashed] (27) -- (33);
        \draw[mydashed] (28) -- (34);
        \draw[mydashed] (29) -- (35);
        \draw[mydashed] (30) -- (36);

    \end{tikzpicture}
}
    \caption{Deterministic pipe route.}
    \label{fig:3}
\end{subfigure}
\hfill
\begin{subfigure}{0.45\linewidth}
        \centering
        \resizebox{\linewidth}{!}{
    \begin{tikzpicture}[node distance={1.5cm}, mycircle/.style={draw, circle, fill=gray!20, thick, minimum size=1cm}, mydashed/.style={gray!50, dashed, line width=3}, mydashed2/.style={gray, dashed, line width=3, black}, mysolid/.style={black, line width=3}, singlewalled1/.style={purple, line width=3}, doublewalled1/.style={green, line width=3}, doublewalled2/.style={green, dashed, line width=3}] 
        \node[mycircle] (1) {1}; 
        \node[mycircle] (2) [right of=1] {2}; 
        \node[mycircle] (3) [right of=2] {3}; 
        \node[mycircle] (4) [right of=3] {4}; 
        \node[mycircle] (5) [right of=4] {5}; 
        \node[mycircle] (6) [right of=5] {6};
        
        \node[mycircle] (7) [below of=1] {7}; 
        \node[mycircle] (8) [right of=7] {8}; 
        \node[mycircle] (9) [right of=8] {9}; 
        \node[mycircle] (10) [right of=9] {10}; 
        \node[] (11) [right of=10] {}; 
        \node[mycircle] (12) [right of=11] {12};
        
        \node[mycircle] (13) [below of=7] {13}; 
        \node[mycircle] (14) [right of=13] {14}; 
        \node[] (15) [right of=14] {}; 
        \node[mycircle] (16) [right of=15] {16}; 
        \node[mycircle] (17) [right of=16] {17}; 
        \node[mycircle] (18) [right of=17] {18};
        
        \node[mycircle] (19) [below of=13] {19}; 
        \node[mycircle] (20) [right of=19] {20}; 
        \node[] (21) [right of=20] {}; 
        \node[mycircle] (22) [right of=21] {22}; 
        \node[mycircle] (23) [right of=22] {23}; 
        \node[mycircle] (24) [right of=23] {24};
        
        \node[mycircle] (25) [below of=19] {25}; 
        \node[mycircle] (26) [right of=25] {26}; 
        \node[mycircle] (27) [right of=26] {27}; 
        \node[mycircle] (28) [right of=27] {28}; 
        \node[mycircle] (29) [right of=28] {29}; 
        \node[mycircle] (30) [right of=29] {30};
        
        \node[mycircle] (31) [below of=25] {31}; 
        \node[mycircle] (32) [right of=31] {32}; 
        \node[mycircle] (33) [right of=32] {33}; 
        \node[mycircle] (34) [right of=33] {34}; 
        \node[mycircle] (35) [right of=34] {35}; 
        \node[mycircle] (36) [right of=35] {36};
        
        \draw[mydashed] (1) -- (2);
        \draw[mydashed] (2) -- (3);
        \draw[mydashed] (3) -- (4);
        \draw[mydashed] (4) -- (5);
        \draw[mydashed] (5) -- (6);
        
        \draw[mydashed] (7) -- (8);
        \draw[mydashed] (8) -- (9);
        \draw[mydashed] (9) -- (10);
        
        \draw[mydashed] (13) -- (14);
        \draw[mydashed] (16) -- (17);
        \draw[mydashed] (17) -- (18);

        \draw[mydashed] (19) -- (20);
        \draw[mydashed] (22) -- (23);
        \draw[mydashed] (23) -- (24);
        
        \draw[mydashed] (25) -- (26);
        \draw[singlewalled1] (26) -- (27);
        \draw[singlewalled1] (27) -- (28);
        \draw[mydashed] (28) -- (29);
        \draw[mydashed] (29) -- (30);

        \draw[mydashed] (31) -- (32);
        \draw[mydashed] (32) -- (33);
        \draw[mydashed] (33) -- (34);
        \draw[mydashed] (34) -- (35);
        \draw[mydashed] (35) -- (36);

        \draw[mydashed] (1) -- (7);
        \draw[mydashed] (2) -- (8);
        \draw[mydashed] (3) -- (9);
        \draw[mydashed] (4) -- (10);
        \draw[mydashed] (6) -- (12);

        \draw[mydashed] (7) -- (13);
        \draw[doublewalled1] (8) -- (14);
        \draw[mydashed] (10) -- (16);
        \draw[mydashed] (12) -- (18);

        \draw[mydashed] (13) -- (19);
        \draw[doublewalled1] (14) -- (20);
        \draw[mydashed] (16) -- (22);
        \draw[mydashed] (17) -- (23);
        \draw[mydashed] (18) -- (24);

        \draw[mydashed] (19) -- (25);
        \draw[doublewalled1] (20) -- (26);
        \draw[singlewalled1] (22) -- (28);
        \draw[mydashed] (23) -- (29);
        \draw[mydashed] (24) -- (30);

        \draw[mydashed] (25) -- (31);
        \draw[doublewalled1] (26) -- (32);
        \draw[mydashed] (27) -- (33);
        \draw[mydashed] (28) -- (34);
        \draw[mydashed] (29) -- (35);
        \draw[mydashed] (30) -- (36);

    \end{tikzpicture}
}
        \caption{Robust pipe route.}
        \label{fig:4}
\end{subfigure}
\caption{Example of a deterministic and robust pipe route connecting fuel tanks (vertex $22$ for diesel and vertex $32$ for methanol) with the engine room (vertex $8$). The purple and green lines represent single and double-walled pipes, respectively. The solid lines denote pipes installed in the first stage, whereas the dashed lines denote pipes installed in the second stage if the methanol scenario takes place. The gray dashed lines represent the remaining edges.} 
\label{fig:example_scenarios}
\end{figure}

To illustrate the importance of considering future fuels, we use the example from Figure~\ref{fig:2}. The first two branches from the scenario tree in Figure~\ref{fig:scenario_tree} are considered in this example: transitioning from diesel to diesel or methanol, respectively. For both fuel types, we need to connect the fuel tanks and the engine room with a pipe route. More specifically, vertices $22$ and $32$ in Figure~\ref{fig:2} represent different fuel tanks, diesel, and methanol, respectively, whereas vertex $8$ denotes the engine room. Figures~\ref{fig:3} and~\ref{fig:4} show a deterministic and robust pipe route for this example, respectively. Because Figure~\ref{fig:3} represents the deterministic solution, the pipe route is only optimal for the first stage, in which the diesel tank gets connected to the engine room. However, it is unknown which fuel type will be used in the future and therefore which scenario will be realized. If diesel persists in the future, the pipe route will not need to change. However, in case methanol occurs in the future, we need to build an extra pipe route in order to connect the methanol tank and the engine room, which is denoted by the green dashed line. Depending on the retrofit costs and the probabilities of the second scenario, a different pipe routing in the first stage can be optimal. For example, Figure~\ref{fig:4} shows a robust solution that does not need any installations in the future. Although this route seems inefficient for the diesel scenario, it prepares the ship for a possible transition to methanol in the future. The only way to make the best pipe routing decisions is to include uncertainty in our mathematical models explicitly.

\subsection{Extension of the deterministic problem formulation to SO and RO}
This section extends the deterministic pipe routing problem formulations to SO and RO formulations. As mentioned in Section~\ref{sec:descriptionSPR}, we need to introduce a set of second-stage scenarios $s \in \mathcal{S}$, where each scenario corresponds to one fuel type. We reuse all previously introduced notations for (DO-U) and (DO-D) but then for each future scenario, which will be indicated by superscript $(s)$ for scenario $s$. For example, $I^{(s)} = (G, \mathcal{P}, \mathcal{T}^{(s)}, P^{(s)}, C^{(s)}, E^{(s)})$ denotes a problem instance in scenario $s$ (in which $\mathcal{T}^{(s)}$ for example denotes all the corresponding terminal groups) and $\textbf{S}^{(s)}$ captures the decisions taken in scenario $s$. The cost of installing a pipe $p \in \mathcal{P}$ at a particular edge $(u,v) \in \mathcal{E}$ is given by $c^{(s)}_{puv} = c_{puv} \cdot \lambda^{(s)}$ where $\lambda^{(s)} > 1$ is the inflation rate. These costs are captured in $C^{(s)} = (c^{(s)}_{puv})_{p \in \mathcal{P}, (u,v) \in \mathcal{E}, s \in \mathcal{S}}$. For an overview of the notation for the new problem formulations, see Table~\ref{tab:notation_overview2}.

\begin{longtable}{@{}p{0.2\textwidth}p{0.8\textwidth}@{}}
\caption{Notation overview for (RO), (SO) and their ILO models, which builds on the notation from Table~\ref{tab:notation_overview1a} and~\ref{tab:notation_overview1b}.}\label{tab:notation_overview2}\\
\toprule
Sets &  \\ 
\midrule
$\mathcal{S}$ & Set of scenario indices, indexed from 1 to S.\\
\midrule
\multicolumn{2}{@{}l@{}}{Decision Variables} \\
\midrule
$d$ & Continuous variable that captures the worst-case retrofit costs for the (RO-U) and (RO-D) models.\\
\bottomrule
\end{longtable}

Using the new notation, we get the following robust optimization problem in which we anticipate the worst-case future scenario:
\begin{subequations} \label{eq:RO_problem}
\begin{align}
    \textbf{\mbox{(RO)}}&& \notag\\
    & \min_{\textbf{S} \in \mathcal{M}(I)} \left( F(I, \emptyset, \textbf{S}) + \max_{s \in \mathcal{S}} \min_{\textbf{S}^{(s)} \in \mathcal{M}(I^{(s)})} \left( F(I^{(s)}, \textbf{S}, \textbf{S}^{(s)}) \right) \right)\\
    & = \min_{\textbf{S} \in \mathcal{M}(I)} \Bigl ( \underbrace{\sum_{(p, (u,v)) \in \textbf{S} \setminus \textbf{S}_0}  c_{puv}}_{\text{First stage costs}} + \underbrace{\max_{s \in \mathcal{S}} \min_{\textbf{S}^{(s)} \in \mathcal{M}(I^{(s)})} \sum_{(p, (u,v)) \in \textbf{S}^{(s)} \setminus \textbf{S}}  c^{(s)}_{puv}}_{\text{Second stage costs}} \Bigr ). \label{eq:RO_problem2}   
\end{align}      
\end{subequations}

We can rewrite (RO) to an undirected ILO model denoted by (RO-U) that (commercial) ILO solvers can solve:
\begin{subequations} \label{eq:basic_robust}
    \begin{alignat}{3}
        \mbox{\textbf{(RO-U)}} && \notag\\
        \min \quad & \sum_{(p, (u,v)) \in (\mathcal{P} \times \mathcal{E})} x_{puv} \cdot c_{puv} + d & \label{eq:basic_ro1} \\
        \mbox{s.t.} \quad & d \geq \sum_{(p, (u,v)) \in (P^{(s)} \times E^{(s)})} (x_{puv}^{(s)} - x_{puv}) \cdot c_{puv}^{(s)} & \qquad \forall s \in \mathcal{S} \label{eq:basic_ro2}\\ 
        \quad & x_{puv}^{(s)} \geq x_{puv} &  \label{eq:basic_ro2}\\
        \quad & \textbf{x}, \textbf{f} \in \mathcal{M}(I)^{\mbox{(DO-U)}} & \label{eq:basic_ro3}\\
        \quad & \textbf{x}^{(s)}, \textbf{f}^{(s)} \in \mathcal{M}(I^{(s)})^{\mbox{(DO-U)}}, & \label{eq:basic_ro4} 
    \end{alignat}
\end{subequations}
where continuous decision variable $d$ captures the worst-case retrofit costs. We force $x_{puv}^{(s)} = 1$ when $x_{puv} = 1$ in~\eqref{eq:basic_ro2}, as we can reuse pipe $p \in \mathcal{P}$ at edge $(u,v) \in \mathcal{E}$ in scenario $s \in \mathcal{S}$. We can also rewrite (RO) to a directed ILO model (RO-D) that is based on (DO-D). So (RO-D) is the same as (RO-U) above, but with $\mathcal{M}(I)^{\mbox{(DO-U)}}$ replaced by $\mathcal{M}(I)^{\mbox{(DO-D)}}$.

The stochastic optimization problem (SO), which aims to minimize the sum of expected costs, can be written as:

\begin{subequations} \label{eq:SO_problem1}
\begin{align}
    \textbf{\mbox{(SO)}}&& \notag\\
    & \min_{\textbf{S} \in \mathcal{M}(I)} \left( F(I, \emptyset, \textbf{S}) + \mathbb{E}_\mathbb{S}\left[ \min_{\textbf{S}' \in \mathcal{M}(I^\mathbb{S})} \left( F(I^\mathbb{S}, \textbf{S}, \textbf{S}') \right) \right] \right)\\
    & = \min_{\textbf{S} \in \mathcal{M}(I)} \Bigl ( \underbrace{\sum_{(p, (u,v)) \in \textbf{S} \setminus \textbf{S}_0}  c_{puv}}_{\text{First stage costs}} + \underbrace{\mathbb{E}_\mathbb{S} \Bigl [ \min_{\textbf{S}' \in \mathcal{M}(I^\mathbb{S})} \sum_{(p, (u,v)) \in \textbf{S}^{(s)} \setminus \textbf{S}}  c^{(s)}_{puv} \Bigr ]}_{\text{Second stage costs}} \Bigr ), \label{eq:SO_problem2}   
\end{align}      
\end{subequations}
where $\mathbb{S}$ denotes a discrete random variable for the future scenario with a known distribution, i.e., $\Pr(\mathbb{S} = s)$ is known for all $s \in \mathcal{S}$.

We can rewrite (SO) to an undirected ILO model denoted by (SO-U) that (commercial) ILO solvers can solve: 
\begin{subequations} \label{eq:basic_stochastic}
    \begin{alignat}{3}
        \mbox{\textbf{(SO-U)}} && \notag\\
        \min \quad & \sum_{(p,(u,v)) \in (\mathcal{P} \times \mathcal{E})} x_{puv} \cdot c_{puv} + \sum_{s \in \mathcal{S}} \rho^{(s)} \sum_{(p, (u,v)) \in (P^{(s)} \times E^{(s)})} (x_{puv}^{(s)} - x_{puv}) \cdot c_{puv}^{(s)} & \label{eq:basic_so1} \\
        \mbox{s.t.} \quad & x_{puv}^{(s)} \geq x_{puv} &  \label{eq:basic_so2}\\
        \quad & \textbf{x}, \textbf{f} \in \mathcal{M}(I)^{\mbox{(DO-U)}} & \label{eq:basic_so3}\\
        \quad & \textbf{x}^{(s)}, \textbf{f}^{(s)} \in \mathcal{M}(I^{(s)})^{\mbox{(DO-U)}}. & \label{eq:basic_so4} 
    \end{alignat}
\end{subequations}
where $\rho^{(s)}$ represents the probability of scenario $s \in \mathcal{S}$. We can also write (SO) to a directed ILO model (SO-D) that is based on (DO-D). So (SO-D) is the same as (SO-U), but with $\mathcal{M}(I)^{\mbox{(DO-U)}}$ replaced by $\mathcal{M}(I)^{\mbox{(DO-D)}}$.

\section{Results} \label{sec:results}
To highlight different aspects of the 2S-SSFP and the proposed models, we consider three graphs of different sizes and with different levels of modeling realism. This section will discuss our experiments on these graphs and report the corresponding results. Section~\ref{sec:smallinstance} elaborates on the small graph from Figure~\ref{fig:2}, whereas Section~\ref{sec:artificialinstance} and Section~\ref{sec:realisticinstance} discuss the artificial and realistic graphs, respectively. We compare the scalability of the models and study the relative gains of taking uncertainty into account. Both (DO-U)~\eqref{eq:basic_deterministic} and (DO-D)~\eqref{eq:advanced_deterministic} are deterministic models that do not consider the uncertainty of the second stage and only optimize for the first stage, which is a naive approach compared to SO and RO. In the following, DO will mainly serve as a benchmark to compare with SO and RO. All experiments are executed on a supercomputer with 32 cores and CPU 2.4GHz using the Gurobi solver \cite{gurobi_optimization_llc_gurobi_2023} for our Python code, which is available upon request.

\subsection{Small graph: an example of the added value from SO} \label{sec:smallinstance}
This section will discuss an example that shows the added value from SO and gives insight into the cost difference between the three mathematical models. We use the small and simple graph from~Figure~\ref{fig:2}, which is a grid graph where the costs of installing a single-walled pipe in the first stage equals $1$, i.e.,
\[c_{1uv} = 1, \qquad \forall (u,v) \in \mathcal{E}.\] 

We assume that double-walled pipes are more expensive than single-walled ones for all three graphs. More specifically, we set the cost ratio between double-walled and single-walled pipes equal to $2$, i.e.,
\[
c_{2uv} = 2c_{1uv}, \qquad \forall (u,v) \in \mathcal{E}.
\]
We assume that the costs of installing pipes are higher in the future than in the present, but we do not know to what extent. Therefore, we use the smallest integer that is greater than $1$ for the ratio between future and present costs of installing pipes, i.e.,
\[
\frac{c^{(s)}_{puv}}{c_{puv}}=\lambda^{(s)}=2, \qquad \forall s \in \mathcal{S}, \forall p \in \mathcal{P}, \forall (u,v) \in \mathcal{E}.
\]
Here, $\lambda^{(s)}=2$ is a constant increase rate for the costs of installing pipes.

\begin{figure}[!h]
    \centering
    \begin{subfigure}{0.40\linewidth}
        \centering
        \resizebox{\linewidth}{!}{
    \begin{tikzpicture}[node distance={1.5cm}, mycircle/.style={draw, circle, fill=gray!20, thick, minimum size=1cm}, mydashed/.style={gray!50, dashed, line width=3}, mydashed2/.style={gray, dashed, line width=3, black}, mysolid/.style={black, line width=3}, singlewalled1/.style={purple, line width=3}, doublewalled1/.style={green, line width=3}, doublewalled2/.style={green, dashed, line width=3}] 
        \node[mycircle] (1) {1}; 
        \node[mycircle] (2) [right of=1] {2}; 
        \node[mycircle] (3) [right of=2] {3}; 
        \node[mycircle] (4) [right of=3] {4}; 
        \node[mycircle] (5) [right of=4] {5}; 
        \node[mycircle] (6) [right of=5] {6};
        
        \node[mycircle] (7) [below of=1] {7}; 
        \node[mycircle] (8) [right of=7] {8}; 
        \node[mycircle] (9) [right of=8] {9}; 
        \node[mycircle] (10) [right of=9] {10}; 
        \node[] (11) [right of=10] {}; 
        \node[mycircle] (12) [right of=11] {12};
        
        \node[mycircle] (13) [below of=7] {13}; 
        \node[mycircle] (14) [right of=13] {14}; 
        \node[] (15) [right of=14] {}; 
        \node[mycircle] (16) [right of=15] {16}; 
        \node[mycircle] (17) [right of=16] {17}; 
        \node[mycircle] (18) [right of=17] {18};
        
        \node[mycircle] (19) [below of=13] {19}; 
        \node[mycircle] (20) [right of=19] {20}; 
        \node[] (21) [right of=20] {}; 
        \node[mycircle] (22) [right of=21] {22}; 
        \node[mycircle] (23) [right of=22] {23}; 
        \node[mycircle] (24) [right of=23] {24};
        
        \node[mycircle] (25) [below of=19] {25}; 
        \node[mycircle] (26) [right of=25] {26}; 
        \node[mycircle] (27) [right of=26] {27}; 
        \node[mycircle] (28) [right of=27] {28}; 
        \node[mycircle] (29) [right of=28] {29}; 
        \node[mycircle] (30) [right of=29] {30};
        
        \node[mycircle] (31) [below of=25] {31}; 
        \node[mycircle] (32) [right of=31] {32}; 
        \node[mycircle] (33) [right of=32] {33}; 
        \node[mycircle] (34) [right of=33] {34}; 
        \node[mycircle] (35) [right of=34] {35}; 
        \node[mycircle] (36) [right of=35] {36};
        
        \draw[mydashed] (1) -- (2);
        \draw[mydashed] (2) -- (3);
        \draw[mydashed] (3) -- (4);
        \draw[mydashed] (4) -- (5);
        \draw[mydashed] (5) -- (6);
        
        \draw[mydashed] (7) -- (8);
        \draw[mydashed] (8) -- (9);
        \draw[mydashed] (9) -- (10);
        
        \draw[mydashed] (13) -- (14);
        \draw[mydashed] (16) -- (17);
        \draw[mydashed] (17) -- (18);

        \draw[mydashed] (19) -- (20);
        \draw[mydashed] (22) -- (23);
        \draw[mydashed] (23) -- (24);
        
        \draw[mydashed] (25) -- (26);
        \draw[singlewalled1] (26) -- (27);
        \draw[singlewalled1] (27) -- (28);
        \draw[mydashed] (28) -- (29);
        \draw[mydashed] (29) -- (30);

        \draw[mydashed] (31) -- (32);
        \draw[mydashed] (32) -- (33);
        \draw[mydashed] (33) -- (34);
        \draw[mydashed] (34) -- (35);
        \draw[mydashed] (35) -- (36);

        \draw[mydashed] (1) -- (7);
        \draw[mydashed] (2) -- (8);
        \draw[mydashed] (3) -- (9);
        \draw[mydashed] (4) -- (10);
        \draw[mydashed] (6) -- (12);

        \draw[mydashed] (7) -- (13);
        \draw[doublewalled1] (8) -- (14);
        \draw[mydashed] (10) -- (16);
        \draw[mydashed] (12) -- (18);

        \draw[mydashed] (13) -- (19);
        \draw[doublewalled1] (14) -- (20);
        \draw[mydashed] (16) -- (22);
        \draw[mydashed] (17) -- (23);
        \draw[mydashed] (18) -- (24);

        \draw[mydashed] (19) -- (25);
        \draw[doublewalled1] (20) -- (26);
        \draw[singlewalled1] (22) -- (28);
        \draw[mydashed] (23) -- (29);
        \draw[mydashed] (24) -- (30);

        \draw[mydashed] (25) -- (31);
        \draw[doublewalled2] (26) -- (32);
        \draw[mydashed] (27) -- (33);
        \draw[mydashed] (28) -- (34);
        \draw[mydashed] (29) -- (35);
        \draw[mydashed] (30) -- (36);

    \end{tikzpicture}
}
        \caption{}
        \label{fig:5}
    \end{subfigure}
    \hfill
    \begin{subfigure}{0.55\linewidth}
    \centering
        \begin{tikzpicture}
        \begin{axis}[
            axis lines = left,
            xlabel = {$\rho^{(2)}$},
            ylabel = {Expected costs},
            ymin=0,
            legend pos=south east,
            width=\linewidth,
        ]
        \addplot [
            domain=0:1, 
            samples=100, 
            color=red,
            line width = 3,
        ]
        {4+16*x};
        \addlegendentry{Route~\ref{fig:3}}
        
        \addplot [
            domain=0:1, 
            samples=100, 
            color=blue,
            line width = 3,
        ]
        {11};
        \addlegendentry{Route~\ref{fig:4}}
        \addplot [
            domain=0:1, 
            samples=100, 
            color=black,
            line width = 3,
        ]
        {9+4*x};
        \addlegendentry{Route~\ref{fig:5}}
        \addplot [
            domain=0:1, 
            samples=100, 
            color=orange,
            line width = 3,
        ]
        {min(4+16*x, 11, 9+4*x)};
        \addlegendentry{SO~\eqref{eq:SP_example}}
        \addplot [
            name path=line1,
            color=black,
            line width = 1,
            dashed
        ]
        coordinates{(0.417, 0)(0.417, 20)};
        \addplot [
            name path=line2,
            color=black,
            line width = 1,
            dashed
        ]
        coordinates{(0.5, 0)(0.5, 20)};
        \addplot [
            thick,
            color=gray,
            fill=gray, 
            fill opacity=0.25
        ]
        fill between[
            of=line1 and line2,
        ];
        \end{axis}
        \end{tikzpicture}    
        \caption{}
        \label{fig:solution_example}
    \end{subfigure}
    \caption{Figure~\ref{fig:5} shows the third pipe route option for the example of Figure~\ref{fig:2}. Note the similarity with Figure~\ref{fig:4}, except the dashed line between vertices $26$ and $32$, which denotes that we possibly install a pipe there in the second stage. Figure~\ref{fig:solution_example} shows the expected costs of the three routes and~\eqref{eq:SP_example} in Figure~\ref{fig:2}. The shaded area indicates the interval of $\rho^{(2)}$ values for which $\mathbb{E} R_{\text{SO}} (p_2) \leq \min \left \{ \mathbb{E} R_{1} (p_2), \mathbb{E} R_{2} (p_2) \right \}$.}
    \label{fig:value_sp}
\end{figure}
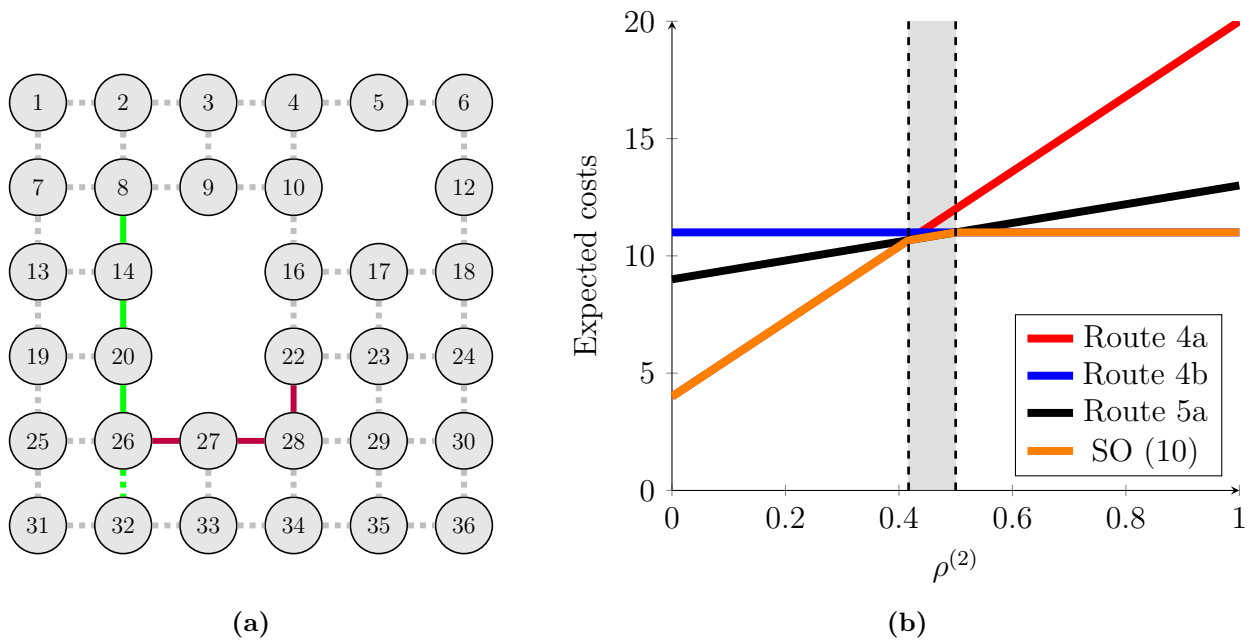

In the example from~Figure~\ref{fig:2}, we start with a diesel scenario in the first stage and transition to either a diesel or methanol scenario in the second stage. In the following, $\rho^{(2)}$ denotes the probability of the methanol scenario in the future, which consists of a transition from diesel to methanol. Since there are only two scenarios in this example, the likelihood of the first scenario is also known: $\rho^{(1)}=1-\rho^{(2)}$. We can algebraically compute the DO, RO, and SO costs for this example. There are three pipe routes for the problem instance in Figure~\ref{fig:2}. The first two routes are DO and RO, as shown in Figures~\ref{fig:3} and~\ref{fig:4}. DO installs four single-walled pipes in the first stage and possibly four double-walled pipes in the second stage, whereas RO installs three single-walled pipes and four double-walled pipes in the first stage. Figure~\ref{fig:5} shows the third route, which installs three single-walled and three double-walled pipes in the first stage and one double-walled pipe in the second stage in the case of the methanol scenario. This route saves costs in the first stage but still prepares for a possible transition in the second stage. As mentioned, the second stage costs of installing a pipe are twice as high compared to the first stage. For a fair comparison, we also consider the expected second-stage costs of installing pipes for DO, which we compute by multiplying the second-stage costs by the probability of the methanol scenario. As a result, the expected DO costs increase linearly with $\rho^{(2)}$. The sum of the expected first and second-stage costs for the three routes are:
\begin{align*}
    & \mathbb{E} R_{1} (\rho^{(2)}) = 4 + 16 \rho^{(2)}\\
    & \mathbb{E} R_{2} (\rho^{(2)}) = 11\\
    & \mathbb{E} R_{3} (\rho^{(2)}) = 9 + 4 \rho^{(2)},
\end{align*}
where $R_{1}$ and $R_{2}$ denote the objective values for the DO and RO solution, and $R_3$ represents the objective value for the solution from Figure~\ref{fig:5}. Note that both $\lambda^{(s)}$ and $\rho^{(2)}$ determine the slope of $R_{1}$ and $R_{3}$. Figure~\ref{fig:solution_example} visualizes the expected costs, where the x-axis and the y-axis show $\rho^{(2)}$ and the sum of the expected first- and second-stage costs, respectively. RO yields a horizontal line because this model does not depend on probabilities.

SO takes the minimum of all three options: 
\begin{equation}
    \mathbb{E} R_{\text{SO}} (\rho^{(2)}) = \min \left \{ 4 + 16 \rho^{(2)}, 11, 9 + 4 \rho^{(2)} \right \}. \label{eq:SP_example}
\end{equation}
$\mathbb{E} R_3$ intersects with $\mathbb{E} R_{1}$ at $\rho^{(2)}=\frac{5}{12}$ and with $\mathbb{E} R_{2}$ at $\rho^{(2)}=\frac{1}{2}$. Between these two values, it holds that $\mathbb{E} R_{\text{SP}} < \min \left \{ \mathbb{E} R_{1}, \mathbb{E} R_{2} \right \}$, which shows the value of taking uncertainty into account. However, this result is based on the assumption that $\rho^{(2)}$ is known, which is not always the case in practice. Another way to quantify the added value of SO is the value of the stochastic solution (VSS), which is the difference between the expectation of the expected value solution (EEVS) and the optimal objective value of the SO. In this example, VSS equals the difference between $\mathbb{E} R_{1}$ and $\mathbb{E} R_{\text{SO}}$, which ranges between $0$ and $9$ ($82\%$ of the SO objective). In other words, the relative gains of SO are relatively high in this example. 

\subsection{Artificial instance: a comparison of the models} \label{sec:artificialinstance}
In this section, we compare the scalability of the three models and show the relative gains of considering uncertainty. To this end, we introduce an artificial graph, which is a 5-by-5 grid graph. The costs of installing a single-walled pipe are randomly drawn from a uniform distribution between $1$ and $10$, i.e.,
\[
c_{1uv} \sim \mathcal{U}[1, 10], \qquad \forall (u,v) \in \mathcal{E}.
\] 
The installation costs of double-walled pipes are again $c_{2uv} = 2c_{1uv}$, $\forall (u,v) \in \mathcal{E}$. The terminals are randomly chosen from all vertices. For a parameter study, we vary the number of scenarios between two and four, the number of terminals between three and five, and the number of terminal groups between one and three, which amounts to $27$ different parameter settings in total for the artificial graph. To account for randomness, we generate $100$ instances for each parameter setting, resulting in $2700$ different instances for the same artificial graph. For the sake of simplicity, we assume that every terminal group can install all pipe types on every edge of the graph and that each scenario is equally probable, i.e., $\rho_s = \frac{1}{S}, \forall s \in \mathcal{S}$.

\begin{figure}[!h]
    \begin{subfigure}{0.40\linewidth}
\begin{tikzpicture}

\definecolor{brown}{RGB}{165,42,42}
\definecolor{darkgray176}{RGB}{176,176,176}
\definecolor{lightblue}{RGB}{173,216,230}

\begin{axis}[
tick align=outside,
tick pos=left,
x grid style={darkgray176},
xlabel={Method},
xmin=-0.2425, xmax=1.2425,
xtick style={color=black},
xtick={0,1},
xticklabels={(SO-U),(SO-D)},
y grid style={darkgray176},
ylabel={Time (sec)},
ymin=0, ymax=6.19728389211496,
ytick style={color=black},
width=\linewidth,
]
\draw[draw=none,fill=lightblue] (axis cs:-0.175,0) rectangle (axis cs:0.175,0.128183985727805);
\addlegendimage{ybar,ybar legend,draw=none,fill=lightblue}
\addlegendentry{Compilation time}

\draw[draw=none,fill=lightblue] (axis cs:0.825,0) rectangle (axis cs:1.175,0.355486509711654);
\draw[draw=none,fill=brown] (axis cs:-0.175,0.128183985727805) rectangle (axis cs:0.175,5.90217513534758);
\addlegendimage{ybar,ybar legend,draw=none,fill=brown}
\addlegendentry{Run time}

\draw[draw=none,fill=brown] (axis cs:0.825,0.355486509711654) rectangle (axis cs:1.175,1.89397803350731);
\end{axis}

\end{tikzpicture}
        \subcaption{}
        \label{fig:runcompilation}
    \end{subfigure}
    \begin{subfigure}{0.40\linewidth}
        \input{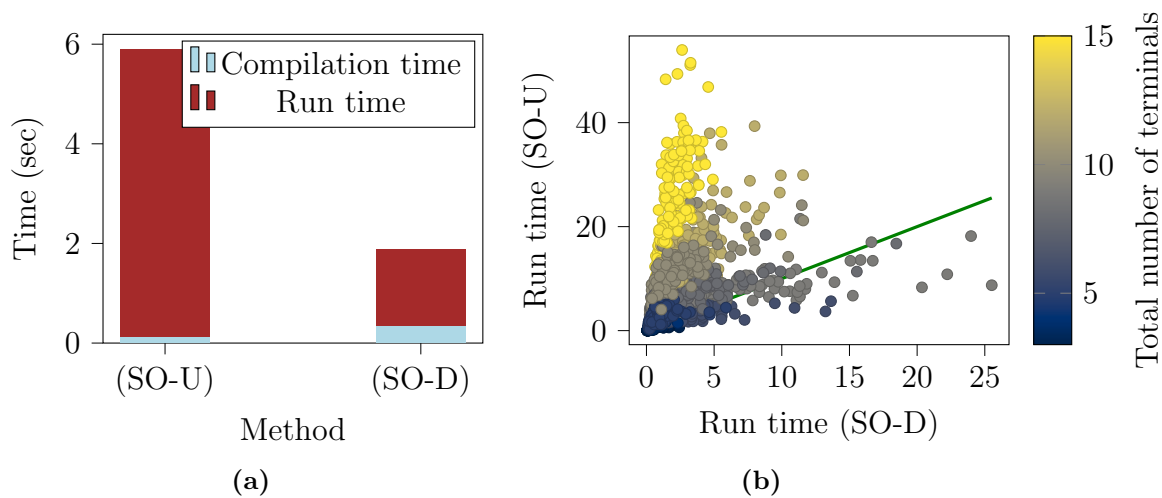}
        \subcaption{}
        \label{fig:scatter}        
    \end{subfigure}
    \caption{Figure~\ref{fig:runcompilation} shows the average of the compilation and run times from (SO-U) and (SO-D). Figure~\ref{fig:scatter} shows the run times of (SO-U) and (SO-D) with the total number of terminals per scenario (i.e., the number of terminal groups multiplied by the number of terminals per terminal group). The green line indicates where the run times of (SO-U) and (SO-D) are equal.}
    \label{fig:runtimes}
\end{figure}

Figure~\ref{fig:runcompilation} shows the average compilation and run time for (SO-U) and (SO-D) over all the $2700$ instances, respectively. The figure shows that directed formulations require more compilation time but yield a considerably shorter run time. We zoom in on this statement with Figure~\ref{fig:scatter}, which displays the run times of (SO-U) and (SO-D) with the total number of terminals per scenario (i.e., the number of terminal groups multiplied by the number of terminals per terminal group). We see that (SO-D) is generally faster than (SO-U), typically when the number of terminals is large.

It is interesting to zoom in on the worst- and best-case instances in terms of run time for (SO-U). The best-case instance contains two scenarios, one terminal group per scenario, and three terminals per terminal group. These terminals lie close to each other, resulting in shorter run times. The worst-case instance contains more terminals, terminal groups, and scenarios, which leads to longer run times, probably due to an increased number of decision variables and constraints. When the terminals are spread across the whole graph, the routing problem becomes more complex, resulting in longer run times.

Because DO, RO, and SO are entirely different models, comparing them with one measure is not trivial. Therefore, we compare them to each other's objective in Table~\ref{tab:comparison_table}. In the DO objective, we divide the first-stage costs of the three models by the first-stage costs of DO. We see that RO yields a higher ratio than SO, which is caused by the conservatism of RO. SO and RO are 29\% and 62\% more expensive than DO in the DO objective, respectively. In the RO objective, we compute the objective value of the three solutions in case of the worst-case scenario and divide it by the RO objective value. We find that the objectives of SO and RO are closer to each other than those of DO and RO. Compared to RO, DO and SO are 29\% and 9\% more expensive, respectively. Hence, considering the uncertainty for the worst case yields a relative gain of 22\%. Finally, in the SO objective, we compute the objective value of the three solutions in the average case with equal probabilities for each scenario and divide it by the SO objective value. DO and RO perform comparably in this objective as both lie approximately 5\% from the optimal objective value. In other words, the VSS amounts to approximately 5\%, which is the expected gain from solving the SO instead of DO.

\begin{table}[!h]
    \centering
    \begin{tabular}{@{}lrrr@{}}
\toprule
Model \textbackslash \; objective &  DO &  RO &  SO \\
\midrule
DO     &      1.000 &      1.286 &      1.057 \\
RO     &      1.623 &      1.000 &      1.048 \\
SO     &      1.225 &      1.116 &      1.000 \\
\bottomrule
\end{tabular}

    \caption{Comparison of the models' average relative performances for three different objectives in artificial instances.}
    \label{tab:comparison_table}
\end{table}

\begin{figure}[!h]
    \centering
    \begin{subfigure}{0.45\linewidth}
    \centering
\begin{tikzpicture}

\definecolor{darkgray176}{RGB}{176,176,176}

\begin{axis}[
tick align=outside,
tick pos=left,
x grid style={darkgray176},
xlabel={Objective RO/objective DO},
xmin=0.625004089325536, xmax=8.87491412416373,
xtick style={color=black},
y grid style={darkgray176},
ylabel={Frequency},
ymin=0, ymax=633.15,
ytick style={color=black},
width=\linewidth,
]
\draw[draw=none,fill=black] (axis cs:1,0) rectangle (axis cs:1.14999836426979,149);
\draw[draw=none,fill=black] (axis cs:1.14999836426979,0) rectangle (axis cs:1.29999672853957,476);
\draw[draw=none,fill=black] (axis cs:1.29999672853957,0) rectangle (axis cs:1.44999509280936,603);
\draw[draw=none,fill=black] (axis cs:1.44999509280936,0) rectangle (axis cs:1.59999345707914,442);
\draw[draw=none,fill=black] (axis cs:1.59999345707914,0) rectangle (axis cs:1.74999182134893,316);
\draw[draw=none,fill=black] (axis cs:1.74999182134893,0) rectangle (axis cs:1.89999018561871,211);
\draw[draw=none,fill=black] (axis cs:1.89999018561871,0) rectangle (axis cs:2.0499885498885,149);
\draw[draw=none,fill=black] (axis cs:2.0499885498885,0) rectangle (axis cs:2.19998691415828,91);
\draw[draw=none,fill=black] (axis cs:2.19998691415828,0) rectangle (axis cs:2.34998527842807,66);
\draw[draw=none,fill=black] (axis cs:2.34998527842807,0) rectangle (axis cs:2.49998364269785,54);
\draw[draw=none,fill=black] (axis cs:2.49998364269785,0) rectangle (axis cs:2.64998200696764,31);
\draw[draw=none,fill=black] (axis cs:2.64998200696764,0) rectangle (axis cs:2.79998037123742,30);
\draw[draw=none,fill=black] (axis cs:2.79998037123742,0) rectangle (axis cs:2.94997873550721,18);
\draw[draw=none,fill=black] (axis cs:2.94997873550721,0) rectangle (axis cs:3.099977099777,9);
\draw[draw=none,fill=black] (axis cs:3.099977099777,0) rectangle (axis cs:3.24997546404678,8);
\draw[draw=none,fill=black] (axis cs:3.24997546404678,0) rectangle (axis cs:3.39997382831657,8);
\draw[draw=none,fill=black] (axis cs:3.39997382831657,0) rectangle (axis cs:3.54997219258635,13);
\draw[draw=none,fill=black] (axis cs:3.54997219258635,0) rectangle (axis cs:3.69997055685614,3);
\draw[draw=none,fill=black] (axis cs:3.69997055685614,0) rectangle (axis cs:3.84996892112592,3);
\draw[draw=none,fill=black] (axis cs:3.84996892112592,0) rectangle (axis cs:3.99996728539571,6);
\draw[draw=none,fill=black] (axis cs:3.99996728539571,0) rectangle (axis cs:4.14996564966549,4);
\draw[draw=none,fill=black] (axis cs:4.14996564966549,0) rectangle (axis cs:4.29996401393528,2);
\draw[draw=none,fill=black] (axis cs:4.29996401393528,0) rectangle (axis cs:4.44996237820506,1);
\draw[draw=none,fill=black] (axis cs:4.44996237820506,0) rectangle (axis cs:4.59996074247485,1);
\draw[draw=none,fill=black] (axis cs:4.59996074247485,0) rectangle (axis cs:4.74995910674463,1);
\draw[draw=none,fill=black] (axis cs:4.74995910674463,0) rectangle (axis cs:4.89995747101442,1);
\draw[draw=none,fill=black] (axis cs:4.89995747101442,0) rectangle (axis cs:5.0499558352842,1);
\draw[draw=none,fill=black] (axis cs:5.0499558352842,0) rectangle (axis cs:5.19995419955399,1);
\draw[draw=none,fill=black] (axis cs:5.19995419955399,0) rectangle (axis cs:5.34995256382378,0);
\draw[draw=none,fill=black] (axis cs:5.34995256382378,0) rectangle (axis cs:5.49995092809356,0);
\draw[draw=none,fill=black] (axis cs:5.49995092809356,0) rectangle (axis cs:5.64994929236335,0);
\draw[draw=none,fill=black] (axis cs:5.64994929236335,0) rectangle (axis cs:5.79994765663313,0);
\draw[draw=none,fill=black] (axis cs:5.79994765663313,0) rectangle (axis cs:5.94994602090292,0);
\draw[draw=none,fill=black] (axis cs:5.94994602090292,0) rectangle (axis cs:6.0999443851727,0);
\draw[draw=none,fill=black] (axis cs:6.0999443851727,0) rectangle (axis cs:6.24994274944249,0);
\draw[draw=none,fill=black] (axis cs:6.24994274944249,0) rectangle (axis cs:6.39994111371227,0);
\draw[draw=none,fill=black] (axis cs:6.39994111371227,0) rectangle (axis cs:6.54993947798206,1);
\draw[draw=none,fill=black] (axis cs:6.54993947798206,0) rectangle (axis cs:6.69993784225184,0);
\draw[draw=none,fill=black] (axis cs:6.69993784225184,0) rectangle (axis cs:6.84993620652163,0);
\draw[draw=none,fill=black] (axis cs:6.84993620652163,0) rectangle (axis cs:6.99993457079141,0);
\draw[draw=none,fill=black] (axis cs:6.99993457079141,0) rectangle (axis cs:7.1499329350612,0);
\draw[draw=none,fill=black] (axis cs:7.1499329350612,0) rectangle (axis cs:7.29993129933098,0);
\draw[draw=none,fill=black] (axis cs:7.29993129933099,0) rectangle (axis cs:7.44992966360077,0);
\draw[draw=none,fill=black] (axis cs:7.44992966360077,0) rectangle (axis cs:7.59992802787056,0);
\draw[draw=none,fill=black] (axis cs:7.59992802787056,0) rectangle (axis cs:7.74992639214034,0);
\draw[draw=none,fill=black] (axis cs:7.74992639214034,0) rectangle (axis cs:7.89992475641013,0);
\draw[draw=none,fill=black] (axis cs:7.89992475641013,0) rectangle (axis cs:8.04992312067991,0);
\draw[draw=none,fill=black] (axis cs:8.04992312067991,0) rectangle (axis cs:8.1999214849497,0);
\draw[draw=none,fill=black] (axis cs:8.1999214849497,0) rectangle (axis cs:8.34991984921948,0);
\draw[draw=none,fill=black] (axis cs:8.34991984921948,0) rectangle (axis cs:8.49991821348927,1);
\end{axis}

\end{tikzpicture}
        \subcaption{}
        \label{fig:ROvsDOhist}
    \end{subfigure}
    \hfill
    \begin{subfigure}{0.45\linewidth}
    \centering
        \input{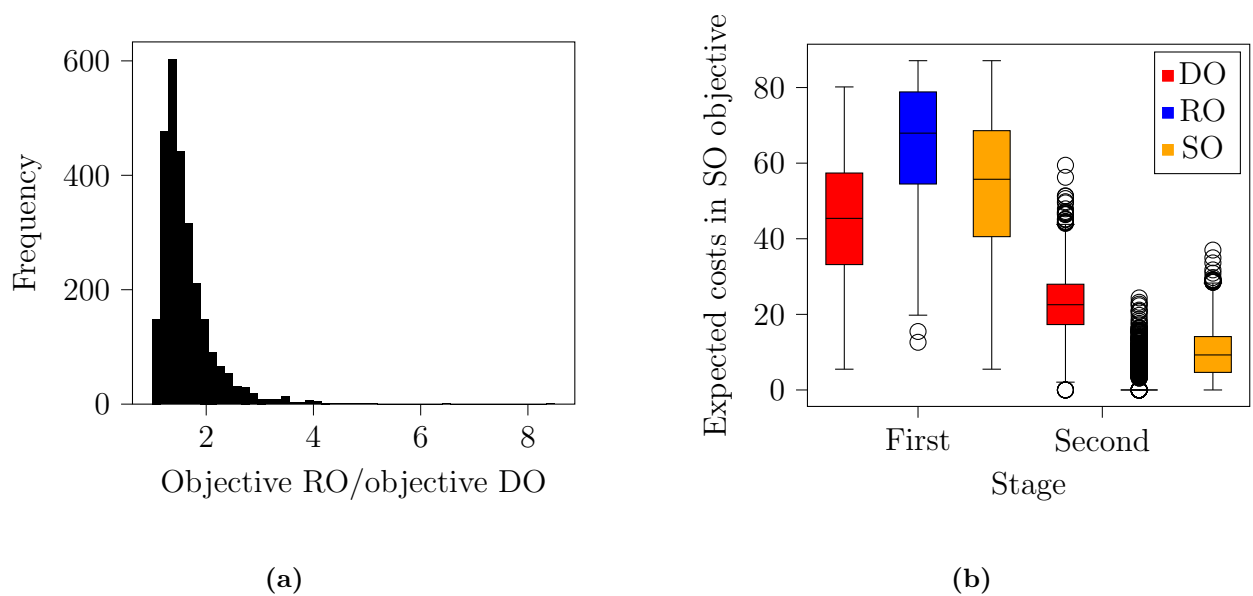}
        \subcaption{}
        \label{fig:SPboxplot}
    \end{subfigure}
    \caption{Figure~\ref{fig:ROvsDOhist} shows a distribution of the RO costs divided by the DO costs for all $2700$ instances. Figure~\ref{fig:SPboxplot} shows boxplots that represent a comparison of all three models in the SO objective.}
    \label{fig:comparison_models}
\end{figure}

The following elaborates on the first- and second-stage costs of all three models. We compare the first-stage costs of DO and RO in Figure~\ref{fig:ROvsDOhist}. We divide the RO costs by the DO costs and find a right-skewed distribution, where the lower bound of this ratio yields $1$. Figure~\ref{fig:SPboxplot} shows how the models perform in the first and second stages. Note that we compute the expected costs in the SO objective, which means that we compute the second stage costs for DO and RO in case of equal probabilities for each scenario, i.e., $\rho^{(1)}=\rho^{(2)}=0.5$. We see that DO has the lowest first-stage costs but is relatively expensive in the second stage, whereas both SO and RO require a bigger investment in the first stage but have considerably lower second-stage costs. When the terminals are close to each other in the first stage but far apart in the second stage, DO seems affordable in the first stage and yields high costs in the second stage.

\subsection{Realistic instance: an example of an application} \label{sec:realisticinstance}
\begin{figure}
    \centering
    \begin{minipage}{\linewidth}
    \begin{subfigure}{\linewidth}
        \centering
        \includegraphics[width=\linewidth]{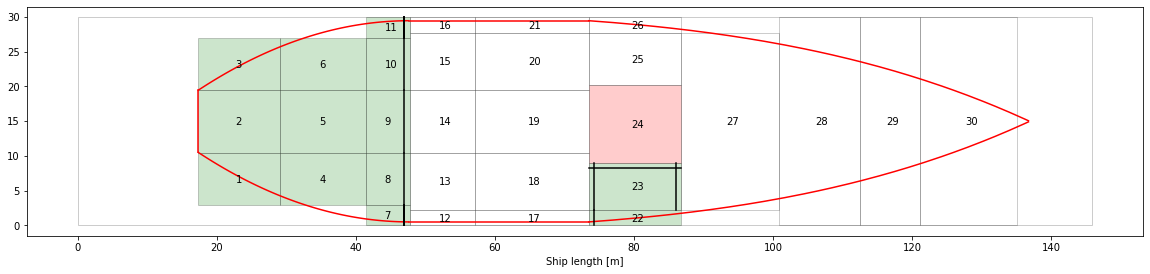}
        \subcaption{Deck 1}
        \label{fig:deck1}
    \end{subfigure}
    \vfill
    \begin{subfigure}{\linewidth}
        \centering
        \includegraphics[width=\linewidth]{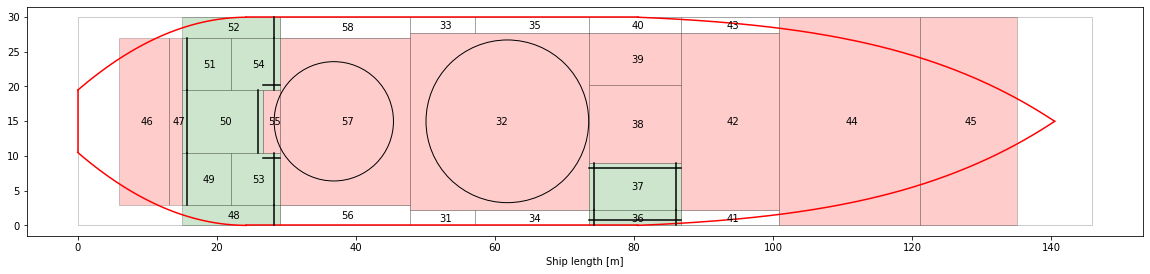}
        \subcaption{Deck 2}
        \label{fig:deck2}
    \end{subfigure}
    \vfill
    \begin{subfigure}{\linewidth}
        \centering
        \includegraphics[width=\linewidth]{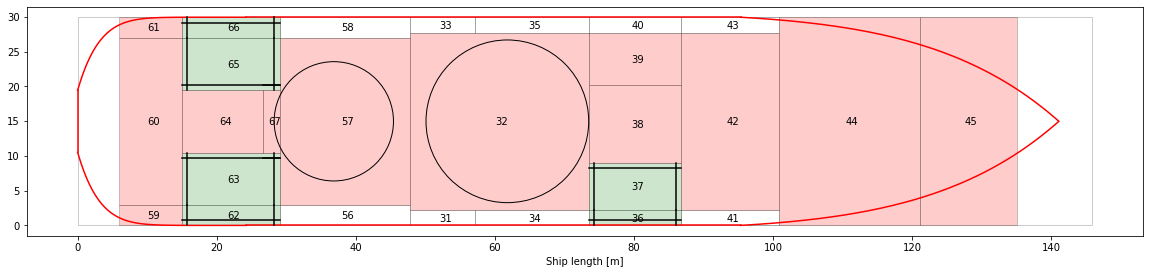}
        \subcaption{Deck 3}
        \label{fig:deck3}
    \end{subfigure}
    \vfill
    \begin{subfigure}{\linewidth}
        \centering
        \includegraphics[width=\linewidth]{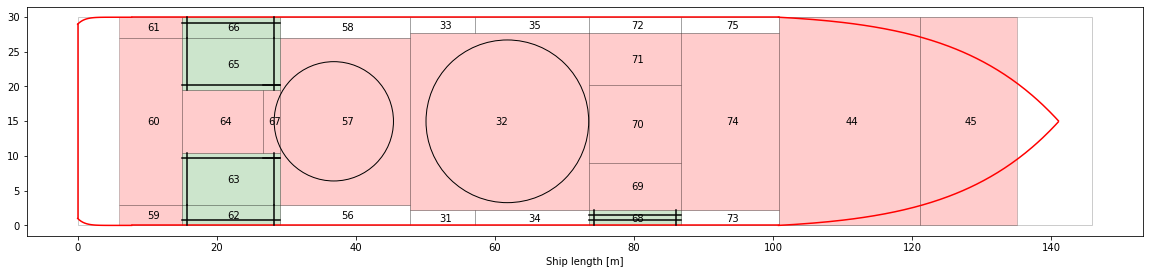}
        \subcaption{Deck 4}
        \label{fig:deck4}
    \end{subfigure}        
    \end{minipage}
    \caption{Figures~\ref{fig:deck1}-\ref{fig:deck4} shows an overview of the four decks of the ship from~\cite{minderhoud_real_2023}, where each number represents a room. Green rooms denote methanol tanks, while red rooms cannot serve as methanol tanks as they are already occupied. Rooms $24$, $38$, and $70$ denote the moonpool.}
    \label{fig:decks}
\end{figure}

In this section, we apply the three models to a realistic graph to show what robust pipe routing looks like in practice. We collaborated with a shipyard to study a ship consisting of four decks. This graph is based on a schematic overview of~\cite{minderhoud_real_2023}, as shown in Figure~\ref{fig:decks}, which computes possible locations of the methanol tanks in a ship currently fueled by diesel. Note that the figures differ slightly from the original paper due to an improvement in the methodology. This ship contains a moonpool, which is an opening in the floor that gives access to the water below, enabling operators to lower tools into the sea. Pipes cannot go through these rooms, making pipe routing more difficult, as it restricts routes. A cargo ship hold would do the same, as pipes cannot enter the cargo space. However, a work ship like this gives more variation since cargo ships have limited spaces and equipment in front of the cargo hold. We represent the 3D network of compartments in the ship as a graph; vertices denote rooms, and edges represent connections between adjacent rooms. The resulting graph contains $75$ vertices and $156$ edges. We use the Manhattan distance between the vertices for the costs $c_{1uv}$ of installing a single-walled pipe in the first stage. The installation costs of double-walled pipes are again $c_{2uv} = 2c_{1uv}$, $\forall (u,v) \in \mathcal{E}$.

We assume that we start with diesel in the first stage and transition to diesel or methanol in the second stage. We use the locations of the current diesel tanks of the ship as terminals for the diesel scenario and the locations of the methanol tanks computed in~\cite{minderhoud_real_2023} as terminals for the methanol scenario. Diesel pipes cannot be routed through the double bottom or through rooms adjacent to the water. Whereas diesel can be routed through either single- or double-walled pipes, methanol requires double-walled pipes due to safety regulations \cite{lloyds_register_rules_2023}. A mathematical overview of the realistic instance is shown below:
\begin{itemize}
    \item $\mathcal{T} = \mathcal{T}^{(1)} = \{37, 42, 53, 54, 63, 65\}$
    \item $\mathcal{T}^{(2)} = \{1-11, 22, 23, 36, 37, 42, 48-54, 62, 63, 65, 66, 68\}$
    \item $\mathcal{P} = P = P^{(1)} = \{1, 2\}$
    \item $P^{(2)} = \{2\}$
    \item $\mathcal{F} = \{1-30, 31, 33, 34, 35, 36, 40, 41, 43, 44-48, 52, 56, 58, 59, 61, 68, 72, 73, 75\}$
    \item $E = E^{(1)} = \left\{ (u,v): u \in \mathcal{V} \setminus \mathcal{F}, v \in \mathcal{V} \setminus \mathcal{F}, u \sim v, u < v \right\}$
    \item $E^{(2)} = \mathcal{E}$,
\end{itemize}
where vertex $42$ represents the engine room and set $\mathcal{F}$ represents the rooms in which we cannot install diesel pipes. We ran all models on the realistic graph and show our findings in Table~\ref{tab:results_RI}. For SO, we assume equal probabilities for the scenarios in the second stage, i.e., $\rho^{(1)}=\rho^{(2)}=\frac{1}{2}$. We see that the models with directed formulations require more compilation time, probably caused by an increase in the number of decision variables and constraints, but yield considerably shorter run times. When focusing on the undirected formulations, we see that the best integer solution is found relatively quickly and that the solver needs a relatively long time to close the gap, whereas directed formulated models close the gap within almost a second.

\begin{table}[!h]
    \centering
    \begin{tabular}{@{}lrrrrrr@{}}
    \toprule
    & \multicolumn{2}{c}{Deterministic} & \multicolumn{2}{c}{Stochastic} & \multicolumn{2}{c}{Robust}\\
    & \multicolumn{1}{r}{(DO-U)} & \multicolumn{1}{r}{(DO-D)} & \multicolumn{1}{r}{(SO-U)} & \multicolumn{1}{r}{(SO-D)} & \multicolumn{1}{r}{(RO-U)} & \multicolumn{1}{r}{(RO-D)}\\ \midrule
    Compilation time & 0.034 & 0.054 & 0.361 & 0.593 & 0.359 & 0.566\\
    Run time & 0.127 & 0.045 & 450.333 & 0.894 & 566.527 & 1.137\\ 
    Decision variables & 1,092 & 1,405 & 10,920 & 12,171 & 10,921 & 12,172\\
    Constraints & 765 & 1,609 & 8,391 & 14,957 & 8,393 & 14,959\\
    \bottomrule
    \end{tabular}
    \caption{Result overview of running all models on the realistic graph.}
    \label{tab:results_RI}
\end{table}

\begin{figure}[!h]
    \centering
    \includegraphics[width=0.45\linewidth]{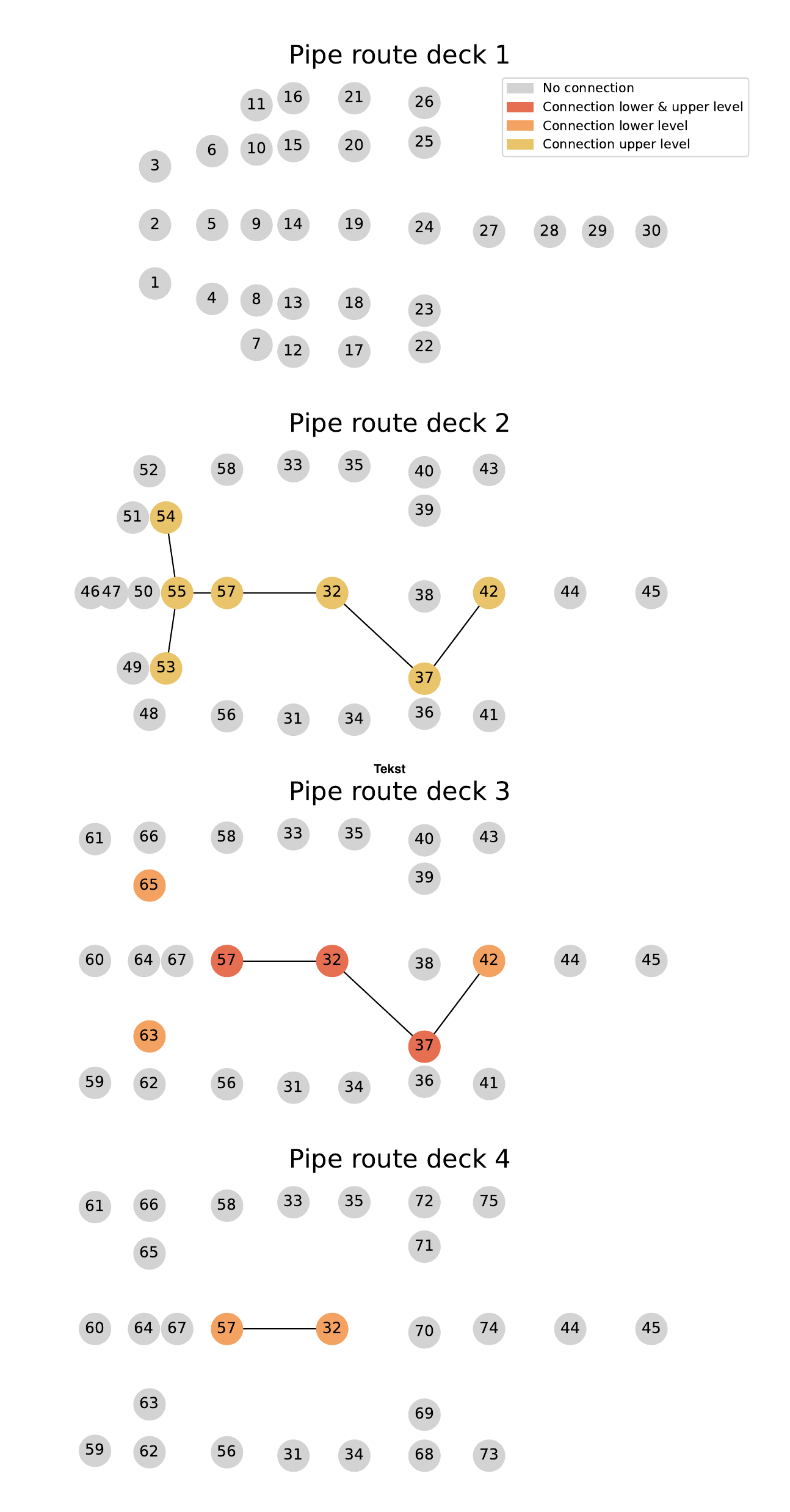}
    \caption{Optimal DO pipe route for the realistic graph, which only consists of single-walled pipes.}
    \label{fig:DO_SW}
\end{figure}

\begin{figure}[!h]
    \centering
    \begin{subfigure}{0.45\linewidth}
        \includegraphics[width=\linewidth]{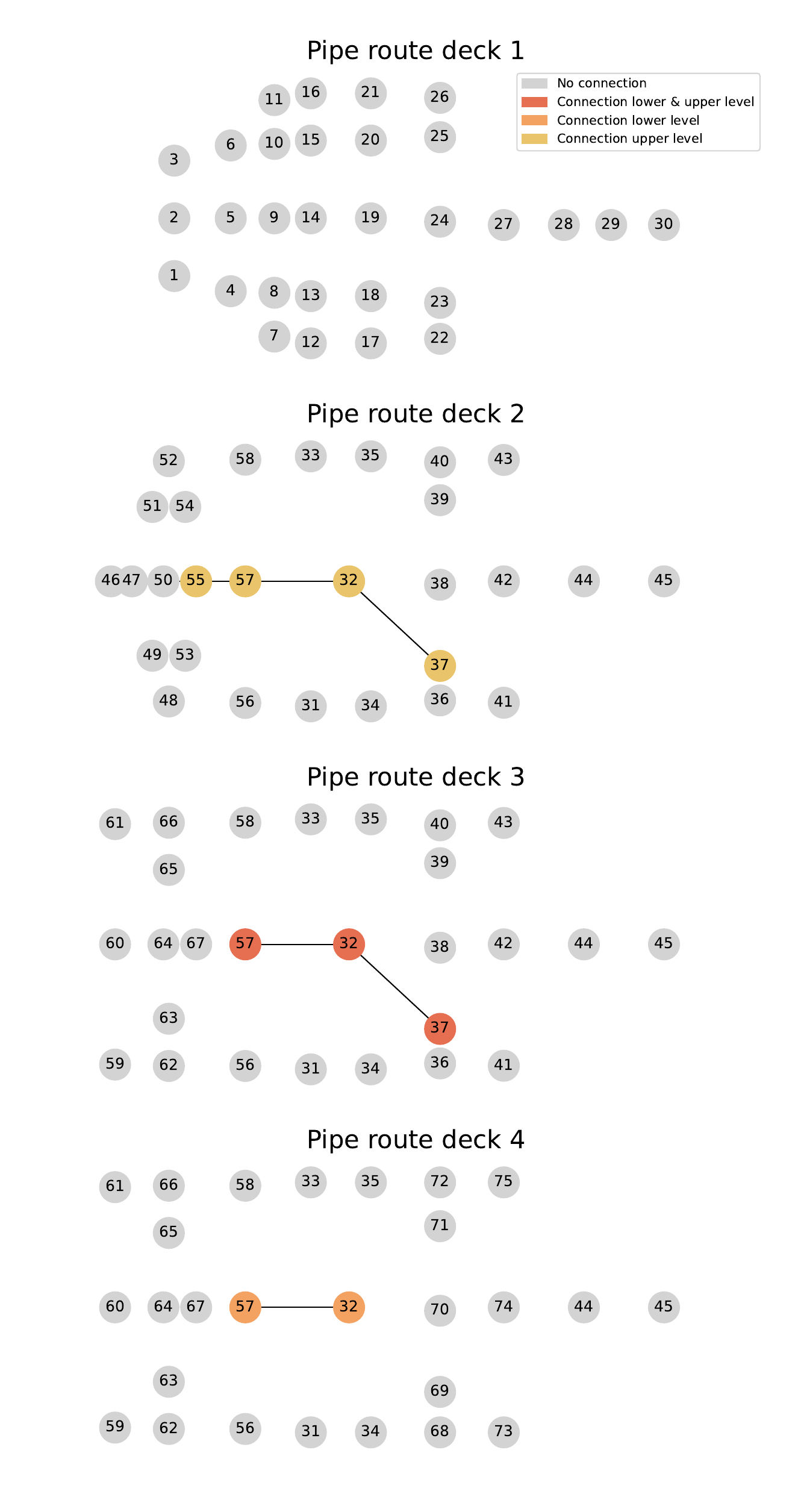}
        \subcaption{}
        \label{fig:RO_SW}
    \end{subfigure}
    \hfill
    \begin{subfigure}{0.45\linewidth}
        \includegraphics[width=\linewidth]{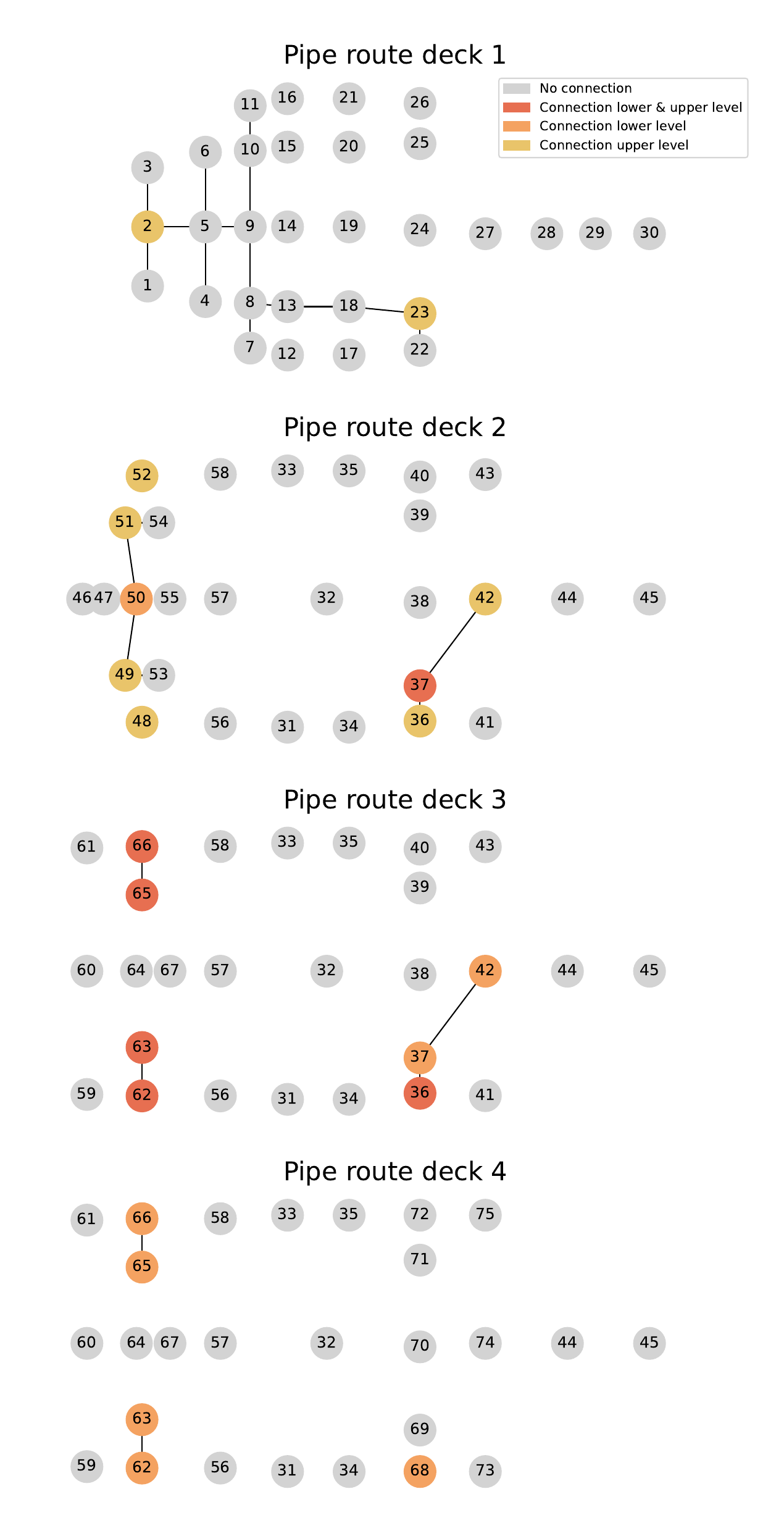}
        \subcaption{}
        \label{fig:RO_DW}
    \end{subfigure}
    \caption{Optimal RO pipe route for the realistic graph. Figure~\ref{fig:RO_SW} shows the single-walled pipe route, whereas Figure~\ref{fig:RO_DW} represents the double-walled pipe route.}
    \label{fig:pipe_route}
\end{figure}

The optimal routes according to DO and RO are displayed in Figure~\ref{fig:pipe_route}. We did not include the optimal route according to SO as it is similar to RO's solution. All three figures use different vertex colors denoting different connections with lower and upper decks. In Figure~\ref{fig:DO_SW}, we see that the DO solution only contains single-walled pipes and avoids rooms adjacent to the water. Figure~\ref{fig:RO_SW} and~\ref{fig:RO_DW} show that the RO solution uses single and double-walled pipes to prepare for diesel and methanol. The RO solution shows three main insights: 1) the pipe route goes via the starboard side of the moonpool to the engine room due to room $37$, which contains a tank as well; 2) the double bottom contains a center pipeline that connects rooms $1$-$11$ with room $23$, which is adjacent to room $37$; and 3) single-walled pipes are placed to connect the diesel tanks to the engine room, as they cannot use the aforementioned pipeline in the double bottom.

\clearpage

\section{Conclusion} \label{sec:conclusion}
Motivated by ship pipe routing under the uncertainty of the energy transition, we have introduced the 2S-SSFP and corresponding DO, SO, and RO models. For each of these three models, we have used both undirected and directed flow formulations based on~\cite{schmidt_stronger_2021}. We have applied them to three graphs to highlight different aspects of the 2S-SSFP and the proposed models.

Our experiments show that considering uncertainty can yield relative gains up to 22\%. An application to a graph representing an actual ship shows three insights: 1) we install double-walled pipes in the center of the double-bottom to connect methanol tanks with the engine room; 2) we use a combination of single- and double-walled pipes on the other decks, and 3) we route the pipe network via the starboard side. Furthermore, we find that DO is typically the cheapest in the first stage but has considerable second-stage costs. SO and RO invest more in the first stage, which results in lower second-stage costs. DO performs the worst when the first stage contains terminals that lie close to each other, whereas the second stage consists of widespread terminals over the graph. The run times of the models increase when the number of scenarios, terminals, and terminal groups increases, especially when the terminals are widespread over the graph. More specifically, the directed formulations require longer compilation times but yield considerably shorter run times than the undirected formulations, which is in line with the findings in \cite{schmidt_stronger_2021}. The difference in run times between directed and undirected formulations increases as the number of terminals increases.

For future research, we suggest studying methods that perform well on larger graphs to take more ship details into account, as an ILO model typically does not scale well. Possible scaleable methods are (meta-)heuristics or the L-shaped method for SO. Additionally, we could make the model more realistic by allowing multiple fuel types within one scenario. The energy transition might consist of multiple stages with different scenarios. Consequently, multi-stage SO could be helpful in this case. Finally, we could make the costs of installing pipes dependent on their location on the ship. For example, installing pipes near an engine room might be more expensive due to safety measures.
\section*{Declaration of interests}
The authors declare that there are no known conflicts of interest that can influence the work presented in this research.

\section*{Acknowledgements}
We thank Jesper Zwaginga for providing the data and helping with the first version of the mathematical models. We thank Joris Slootweg and Ruurd Buijs for their insightful input and feedback. Additionally, we thank SURF (\url{www.surf.nl}) for the support in using the National Supercomputer Snellius. This publication is part of the project READINESS with project number TWM.BL.019.002 of the research program \textit{Topsector Water \& Maritime: the Blue route} which is partly financed by the Dutch Research Council (NWO).

\bibliographystyle{abbrv} 
\bibliography{references}

\end{document}